\documentclass[pdflatex,sn-mathphys-num]{sn-jnl}

\usepackage{graphicx}%
\usepackage{multirow}%
\usepackage{amsmath,amssymb,amsfonts}%
\usepackage{amsthm}%
\usepackage{mathrsfs}%
\usepackage[title]{appendix}%
\usepackage{xcolor}%
\usepackage{textcomp}%
\usepackage{manyfoot}%
\usepackage{booktabs}%
\usepackage{algorithm}%
\usepackage{algorithmicx}%
\usepackage{algpseudocode}%
\usepackage{listings}%
\usepackage[T1]{fontenc}
\usepackage[utf8]{inputenc}

\theoremstyle{thmstyleone}%
\theoremstyle{thmstyletwo}%

\theoremstyle{thmstylethree}%

\newtheorem{theo}{Theorem}
\newtheorem{lem}[theo]{Lemma}

\newtheorem{cor}[theo]{Corollary}
\newtheorem{rem}{Remark}
\newtheorem{exam}{Example}
\newtheorem{definiti}{Definition}

\newenvironment{dem}[1][Proof]{\noindent \textbf{#1.} }{\ \rule{0.5em}{0.5em}}
\begin{document}

\title[Subdifferential of the supremum]{Explicit data-dependent characterizations of the subdifferential of convex pointwise suprema and optimality conditions}


\author[1]{\fnm{Stephanie} \sur{Caro}}\email{stcaro@unap.cl}
\equalcont{These authors contributed equally to this work.}

\author*[3]{\fnm{Abderrahim} \sur{Hantoute}}\email{hantoute@ua.es}
\equalcont{These authors contributed equally to this work.}

\affil[1]{\orgdiv{Faculty of Sciences}, \orgname{Institute of Exact and Natural Sciences, University Arturo Prat}, \orgaddress{
\city{Iquique}, 
\country{Chile}}}

\affil*[2]{\orgdiv{Mathematical Department}, \orgname{University of Alicante}, \orgaddress{
\city{Alicante}, 
\country{Spain}}}



\date{\today}
\maketitle
\begin{abstract}
 
  We establish explicit data-dependent and symmetric characterizations of the
subdifferential of the supremum of convex functions, formulated directly in
terms of the underlying data functions. In our approach, both active and
non-active functions contribute equally through their subdifferentials,
thereby avoiding the need for additional geometric constructions, such as the
domain of the supremum, that arise in previous developments. Applications to
infinite convex optimization yield sharp Karush--Kuhn--Tucker and Fritz--John
optimality conditions, expressed exclusively in terms of the objective and
constraint functions and clearly distinguishing the roles of (almost) active
and non-active constraints.
\end{abstract}

\textbf{keywords:} Supremum of convex functions, subdifferentials, normal
cone to the domain, optimality conditions


\section{Introduction}
Describing the variational properties and subdifferentials of pointwise
suprema in terms of the underlying data functions is a fundamental problem in
optimization and variational analysis, both from theoretical and applied
perspectives. It arises naturally in classical topics such as duality theory,
minimax principles, and the derivation of optimality conditions for
optimization problems with finitely or infinitely many constraints. The study
of this problem dates back to the early developments in convex analysis, with
pioneering contributions in \cite{Br72, Danskin, IoffeLevin, Ro70, Val69,
Vole1, Za02}, and continues to attract considerable attention in more recent
works, including \cite{HLZ08, Ioffe12, chi, Mor13, Vole2, LoVo10, LoVa11}. We
also refer\ to \cite{CH24a, CHLCompactification} for further recent advances
in this area. A comprehensive treatment linking fundamental results of\ convex
analysis to the variational behavior of supremum functions can be found in
\cite{CHLBook}.

The supremum operation has no true counterpart in classical differential
calculus: the pointwise supremum of even two differentiable functions (e.g.,
the absolute value on $\mathbb{R}$) typically fails to be differentiable. In
nonsmooth analysis, however, it can be rigorously treated using generalized
differentiation, including the Fenchel subdifferential in the convex setting
as in \cite{HLZ08, Ioffe12}, and more general nonsmooth subdifferentials such
as those considered in \cite{Morbook, Mor13, PPnoconvex}. 

The subdifferential of a supremum is largely determined by the data functions
that are active---or nearly active---at the reference point. Early results,
such as \cite{IoffeLevin, Val69, Vole1}, involve only these almost active
functions due to continuity assumptions on the supremum. When these
assumptions are relaxed, non-active functions also contribute, typically
through the normal cone to (finite-dimensional sections of) the effective
domain of the supremum (e.g., \cite{HLZ08, Ioffe12, chi, LoVo10}), resulting
in asymmetric formulas.

A natural question is whether symmetric characterizations of the
subdifferential are possible, in which both active and non-active functions
contribute equally through their subdifferentials, thereby avoiding the need
for additional geometric devices\ such as the domain of the supremum, which
can be difficult to handle numerically. This problem was studied in
\cite{CompactCase} (partially extended in \cite{HaLo22}) under a
compact-continuous framework, yielding explicit descriptions that
appropriately weight the contributions of non-active functions. The main goal
of this paper is to provide sharp and operationally useful symmetric
characterizations of the subdifferential, expressed directly in terms of the
underlying data functions, thereby extending the results of \cite{CompactCase}
beyond the compact-continuous framework.

Formally, let $f_{t}:X\rightarrow\overline{{\mathbb{R}}}={\mathbb{R\cup
\{\pm\infty\}}},$ $t\in T,$ be an arbitrary collection\ of lower
semicontinuous convex functions defined on a locally convex space $X.$ We
prove that, for every $x$ in the domain of $f,$ $\operatorname{dom}f$,  the
subdifferential $\partial f$ of the pointwise supremum $f:=\sup_{t\in T}f_{t}$
is completely and exclusively characterized in terms of the\ $\varepsilon
$-subdifferentials\ $\partial_{\varepsilon}f_{t}\ $of the data functions as
follows (Theorem \ref{thmsub}):
\[
\partial f(x)=\bigcap_{\varepsilon>0}\overline{\operatorname*{co}}\left(
A_{\varepsilon}(x)+\varepsilon(B_{\varepsilon,\alpha}(x)\cup C_{\varepsilon
}(x)\cup\{\theta\})\right)  ,
\]
where 

\begin{itemize}
\item $A_{\varepsilon}(x):=%
{\textstyle\bigcup\limits_{t\in T_{\varepsilon}(x)}}
\partial_{\varepsilon}f_{t}(x)$ represents the $\varepsilon$-active data
functions $f_{t}$ at $x,$ that is, those with indices in $T_{\varepsilon
}(x):=\{t\in T:f_{t}(x)\geq f(x)-\varepsilon\};$
\\
\item $B_{\varepsilon,\alpha}(x):=%
{\textstyle\bigcup\limits_{t\in\mathcal{P}\backslash T_{\varepsilon}(x)}}
\partial_{\varepsilon}(\alpha_{t}f_{t})(x),$ with $\mathcal{P}:=\{t\in
T:f_{t}$ is proper$\},$ represents the contribution of the non-$\varepsilon
$-active data functions $f_{t}$ at $x,$ which are penalized by the factor
$\varepsilon>0$ in the above formula;
\\
\item $C_{\varepsilon}(x):=%
{\textstyle\bigcup\limits_{t\in T\backslash\mathcal{P}}}
\mathrm{N}_{\operatorname{dom}f_{t}}^{\varepsilon}(x)$ accounts for\ nonproper
functions having empty $\varepsilon$-subdifferentials but which may still
contribute to $\partial f(x).$\\
\end{itemize}
The parameters $\alpha_{t}>0$ above are suitably chosen weights (see Theorem
\ref{thmsub}). The term $C_{\varepsilon}(x)$ vanishes when all functions
$f_{t}$ are proper. The above representation subsumes most
previously known\ formulas for $\partial f(x)$. For instance, it covers the
general formula given in \cite[Theorem 4]{HLZ08} under the condition
$\operatorname{cl}f=\sup_{t\in T}\operatorname{cl}f_{t},$ relating the closed
hulls of $f$ and the $f_{t}$,
\begin{equation}
\partial f(x)=\bigcap_{L\in\mathcal{F}(x),\ \varepsilon>0}\overline
{\operatorname*{co}}\left(  \bigcup_{t\in T_{\varepsilon}(x)}\partial
_{\varepsilon}f_{t}(x)+\mathrm{N}_{L\cap\operatorname{dom}f}(x)\right)
,\label{refdem16}%
\end{equation}
where $\mathcal{F}(x)$ denotes the family of finite-dimensional linear
subspaces containing $x$. In fact, our approach builds upon this last
characterization, leading to a new representation of the normal cone to
$\operatorname{dom}f.$ Applications of the first formula to infinite convex
optimization give rise to sharp Karush--Kuhn--Tucker and Fritz--John
optimality conditions, expressed exclusively in terms of the objective and
constraint functions and clearly distinguishing the roles of $\varepsilon
$-active and non-$\varepsilon$-active constraints.

This paper is organized as follows. Section \ref{sepre} introduces the
notation and presents the preliminary results used throughout the paper.
Section \ref{NormalCone} is devoted to characterizing the normal cone to the
effective domain $\operatorname{dom}f$; the main result is stated in Theorem
\ref{Theo13}. The desired characterization of $\partial f(x)$ is established
in Theorem \ref{thmsub} in Section \ref{SectionSubdifferential}. Section
\ref{secopt} derives optimality conditions for infinite convex optimization
based on the previously obtained characterizations. Finally, Section
\ref{secultim2} presents variants of the main characterizations of $\partial
f(x)$ under additional continuity assumptions.

\section{Notation and preliminaries\label{sepre}}

We consider\ a real (separated) locally convex space $X$ (lcs, for short),
whose topological dual is denoted by $X^{\ast}.$ The topology in $X^{\ast}$
can be any topology compatible with the duality pair $(X,X^{\ast})$,
associated with\ the bilinear form $(x^{\ast},x)\in X^{\ast}\times
X\mapsto\langle x^{\ast},x\rangle:=x^{\ast}(x).$ Typical examples of such
topologies include the weak* topology ($w^{\ast}$), and the norm topology when
$X$ is a reflexive Banach space. In both $X$\ and $X^{\ast},$ the zero vector
is denoted by $\theta.$ The family of closed, convex, and balanced
neighborhoods of\textbf{\ }$\theta\ $(called\textbf{\ }$\theta$\textbf{-}%
neighborhoods\textbf{)} is denoted\ by\textbf{\ }$\mathcal{N}$.

We use the notation $\overline{\mathbb{R}}:=\mathbb{R}\cup\{-\infty,+\infty
\}$, $\mathbb{R}_{\infty}:=\mathbb{R}\cup\{+\infty\}$, and adopt the
conventions\emph{\ }$\left(  +\infty\right)  +(-\infty)=+\infty,$
$0.(+\infty):=+\infty$, $0.(-\infty):=0$, $\sum\nolimits_{t\in\emptyset}=0,$
$\cup_{\emptyset}\emptyset=\emptyset$. Let\ $\mathbb{R}_{+}^{(T)}%
:=\{\lambda:T\rightarrow\mathbb{R}_{+}:\operatorname*{supp}\!\lambda$
finite$\},$ where $\operatorname*{supp}\!\lambda:=\{t\in T:\lambda(t)\neq0\},$
and define $\Delta(T):=\{\lambda\in\mathbb{R}_{+}^{(T)}:%
{\textstyle\sum_{t\in\operatorname*{supp}\!\lambda}}
\lambda(t)=1\}.$ In particular, the $m$-canonical simplex is $\Delta
_{m}:=\Delta(\{1,\cdots,m\}),$ $m\geq1.$

Given $A$, $B\subset X,$ their Minkowski sum is $A+B:=\{a+b:\ a\in A,\text{
}b\in B\}$, with the convention $A+\emptyset=\emptyset+A=\emptyset$. For
$\Lambda\subset\mathbb{R}$, we write\ $\Lambda A:=\left\{  \lambda
a:\lambda\in\Lambda,\text{ }a\in A\right\}  $ (with $\Lambda\emptyset
=\emptyset A=\emptyset)$. The convex hull, conical convex hull, span,
interior, closed hull, negative dual cone, and orthogonal subspace of $A$ are
denoted by $\operatorname*{co}(A)$, $\operatorname*{cone}(A)$,
$\operatorname*{span}(A)$, $\operatorname{int}(A)$, $\operatorname{cl}(A)$ (or
$\overline{A}$), $A^{-}$, and $A^{\perp}$, respectively. By the Bipolar
Theorem we have $A^{--}:=(A^{-})^{-}=\overline{\operatorname{cone}}(A).$
The\ recession cone of\ a nonempty closed convex set $A$\ is $[A]_{\infty
}:=\cap_{t>0}t(A-a)$ for an arbitrary $a\in A$.$\ $By convention, we set
$[\emptyset]_{\infty}=\{\theta\}.$ Consequently, for any family of closed
convex sets $A_{i}\subset X,$ $i\in I,$ we have $[A_{i}]_{\infty}%
=[\overline{\operatorname*{co}}(A_{i}\cup\{\theta\})]_{\infty}\subset
\overline{\operatorname*{co}}(A_{i}\cup\{\theta\})$ and, hence,
\begin{equation}
\cap_{i\in I}[A_{i}]_{\infty}=\cap_{i\in I}[A_{i}\cup\{\theta\}]_{\infty
}=[\cap_{i\in I}(A_{i}\cup\{\theta\})]_{\infty}=[\cap_{i\in I}A_{i}]_{\infty
}.\label{ricone}%
\end{equation}
For every convex set $A$ and every closed linear subspace $L\subset X,$ we
have that
\begin{equation}
\lbrack\operatorname*{cl}(A+L)]_{\infty}=[\operatorname*{cl}(A+L)]_{\infty
}+L.\label{su}%
\end{equation}
If\ $A\subset X^{\ast}$ is nonempty and $x\in X,$ then
combining\ \cite[Exercise 10]{CHLBook} and \cite[Equation 2.23]{CHLBook}
yields
\begin{equation}%
{\textstyle\bigcap_{L\in\mathcal{F}(x),\text{ }L\subset L_{0}}}
\left[  \overline{\operatorname{co}}(A+L^{\perp})\right]  _{\infty}=\left[
\overline{\operatorname{co}}(A)\right]  _{\infty},\label{lemIntReca}%
\end{equation}
for any $L_{0}\in\mathcal{F}(x):=\{L\subset X:L\text{ is a finite-dimensional
linear subspace with }x\in L\}$.

Given a function $f:X\rightarrow\overline{\mathbb{R}}$, we denote by
$\operatorname*{dom}f$ and $\operatorname*{epi}f$ its (effective) domain and
epigraph,$\ $respectively. We say that $f$ is proper if $\operatorname*{dom}%
f\neq\emptyset$ and $f(x)>-\infty,$ for all $x\in X$. Given $c\in\mathbb{R}$,
we set $[f\leq c]:=\{x\in X:f(x)\leq c\}\ $and $[f<c]:=\{x\in X:f(x)<c\}.$ The
closed hull of\ $f$ is the function $\operatorname*{cl}f$ (or $\overline{f}$)
whose epigraph is $\operatorname*{cl}(\operatorname*{epi}f);$ the function $f$
is lower semicontinuous (lsc, for short) if $f=\bar{f}.$ The class of proper
lsc convex functions on $X$ is denoted by $\Gamma_{0}(X)$. For a set $A\subset
X^{\ast}$ (or $X$), its support function is
\[
\sigma_{A}:=\sup_{a\in A}\left\langle a,\cdot\right\rangle .
\]
The indicator function\ of $A$ is defined by%
\[
\mathrm{I}_{A}(a):=0\text{ if }a\in A,\text{ and }\mathrm{I}_{A}%
(a):=+\infty\text{ if not.}%
\]
We adopt\ the conventions $\sigma_{\emptyset}\equiv-\infty$ and $\mathrm{I}%
_{\emptyset}\equiv+\infty.$

If $A\subset X$ is nonempty, then $\operatorname{cl}(\operatorname{dom}%
\sigma_{A})=([\overline{\operatorname*{co}}\left(  A\right)  ]_{\infty})^{-}.$
If $B\subset X$ is another nonempty set,$\ $we have
\[
\operatorname*{dom}(\max(\sigma_{A},\sigma_{B}))=\operatorname*{dom}%
(\sigma_{A+B})=\operatorname*{dom}\sigma_{A\cup B}=\operatorname*{dom}%
\sigma_{A}\cap\operatorname*{dom}\sigma_{B},
\]
and so\ (see \cite[Exercise 30]{CHLBook})
\begin{equation}
([\overline{\operatorname*{co}}\left(  A+B\right)  ]_{\infty})^{-}%
=([\overline{\operatorname*{co}}\left(  A\cup B\right)  ]_{\infty})^{-}%
\subset\left(  \left[  \overline{\operatorname*{co}}\left(  A\right)  \right]
_{\infty}\right)  ^{-}\cap\left(  \left[  \overline{\operatorname*{co}}\left(
B\right)  \right]  _{\infty}\right)  ^{-}.\label{SupDomNega}%
\end{equation}
More generally, if $C$ is a nonempty bounded set $C\subset X,$ then for every
$\gamma>0$
\begin{equation}
([\overline{\operatorname*{co}}\left(  A\cup B\right)  ]_{\infty}%
)^{-}=([\overline{\operatorname*{co}}\left(  A\cup\gamma B\cup C\right)
]_{\infty})^{-}=([\overline{\operatorname*{co}}\left(  A+\gamma B+C\right)
]_{\infty})^{-}.\label{mostakim}%
\end{equation}
Given $\varepsilon\geq0$, the $\varepsilon$-subdifferential of $f:X\rightarrow
\overline{\mathbb{R}}$ at $x\in X$ is
\[
\partial_{\varepsilon}f(x):=\{x^{\ast}\in X^{\ast}:\langle x^{\ast}%
,y-x\rangle\leq f(y)-f(x)+\varepsilon,\text{ for all }y\in X\},
\]
when $x\in f^{-1}(\mathbb{R}),$ and $\partial_{\varepsilon}f(x):=\emptyset$
otherwise. The (exact) subdifferential of $f$\ at $x$ is $\partial
f(x):=\partial_{0}f(x)=\cap_{\varepsilon>0}\partial_{\varepsilon}f(x).$
Moreover, we have
\begin{equation}
\partial_{\varepsilon}f(x)\subset\partial_{\varepsilon}\bar{f}(x)=\partial
_{\varepsilon+f(x)-\bar{f}(x)}f(x),\label{relsub}%
\end{equation}
and, therefore, the condition $\partial_{\varepsilon}f(x)\neq\emptyset$
implies\ $\bar{f}(x)\geq f(x)-\varepsilon.$ If $f:X\rightarrow\overline
{{\mathbb{R}}}$ is convex with nonempty domain and improper closure, then
(see, e.g., \cite[Lemma 2.20]{CH24a})
\begin{equation}
\operatorname{cl}(f^{+})=(\operatorname{cl}f)^{+}=\mathrm{I}%
_{\operatorname{dom}\left(  \operatorname{cl}f\right)  }%
,\label{PositivePart1a}%
\end{equation}
where $f^{+}:=\max\{f,0\}$ denotes the positive part of $f.$ Consequently, if
$f$ is additionally lsc, then $\operatorname{dom}f=\operatorname{cl}%
(\operatorname{dom}f)$\ and, for every $x\in\operatorname{dom}f$, $u\in X$ and
$\varepsilon>0,$
\begin{equation}
(f^{+})_{\varepsilon}^{\prime}(x;u)=(\mathrm{I}_{\operatorname{dom}%
f})_{\varepsilon}^{\prime}(x;u)=\sigma_{\partial_{\varepsilon}\mathrm{I}%
_{\operatorname{dom}f}(x)}(u),\label{PositivePart1b}%
\end{equation}
where $g_{\varepsilon}^{\prime}(x;u):=\inf_{s>0}\frac{g(x+su)-g(x)+\varepsilon
}{s}$ denotes the $\varepsilon$-directional derivative of a function $g$
defined on $X.$ Moreover, in this latter case, for every $x\in
\operatorname{dom}f,$ $\alpha\in\mathbb{R}$ and $\varepsilon\geq0$ (see
\cite[Corollary 5.1.9]{CHLBook})%
\begin{equation}
\partial_{\varepsilon}(\max\{f,\alpha\})(x)=\partial_{\varepsilon}%
(\mathrm{I}_{\operatorname{dom}f}+\alpha)(x)=\mathrm{N}_{\operatorname{dom}%
f}^{\varepsilon}(x),\label{new}%
\end{equation}
where $\mathrm{N}_{A}^{\varepsilon}(x):=\partial_{\varepsilon}\mathrm{I}%
_{A}(x)$ is\ the $\varepsilon$-normal set to $A\subset X$ at $x\in X$; in
particular, $\mathrm{N}_{A}(x):=\mathrm{N}_{A}^{0}(x)$ is the normal cone\ to
$A\subset X$ at $x$. Note that\ $\mathrm{N}_{A}^{\varepsilon}(x)=\mathrm{N}%
_{\bar{A}}^{\varepsilon}(x)$ whenever $x\in A.$ Consequently,
since\ $\operatorname*{cl}(\operatorname*{dom}f)=\operatorname*{cl}%
(\operatorname*{dom}\bar{f})$ for every function $f$\ on $X$, it follows
that\
\[
\mathrm{N}_{\operatorname*{dom}f}^{\varepsilon}(x)=\mathrm{N}%
_{\operatorname*{dom}\bar{f}}^{\varepsilon}(x),\text{ \ for all }%
x\in\operatorname*{dom}f\text{ and }\varepsilon\geq0.
\]

The following two technical lemmas are used in the sequel.

\begin{lem}
\label{ChangeParameter} Let\ $f_{t}:X\rightarrow\overline{\mathbb{R}}$, $t\in
T,$ be lsc convex functions, and let $A\subset X$ be a convex set.$\ $Fix
$\varepsilon>0.$ Then, for all $\delta>0$ and $\beta_{t}>0$ with
$t\in\mathcal{P}:=\{t\in T:f_{t}$ is proper$\},$
\begin{equation}
A_{\delta}:=\left[  \overline{\operatorname*{co}}\left(  \cup_{0<\alpha
_{t}\leq\delta\beta_{t},\text{ }t\in\mathcal{P}}\partial_{\varepsilon}\left(
\alpha_{t}f_{t}\right)  (x)\cup A\right)  \right]  _{\infty}\subset\left[
\overline{\operatorname*{co}}\left(  \cup_{t\in\mathcal{P}}\partial
_{\varepsilon}\left(  \beta_{t}f_{t}\right)  (x)\cup A\right)  \right]
_{\infty}. \label{one}%
\end{equation}
Consequently, for any scalars $\alpha_{t},\beta_{t}>0$ $(t\in\mathcal{P})$
such that $\sup_{t\in\mathcal{P}}\alpha_{t}\beta_{t}^{-1}<+\infty$ and
$\sup_{t\in\mathcal{P}}\alpha_{t}^{-1}\beta_{t}<+\infty,$ we have\
\begin{equation}
\left[  \overline{\operatorname*{co}}\left(  \cup_{t\in\mathcal{P}}%
\partial_{\varepsilon}\left(  \alpha_{t}f_{t}\right)  (x)\cup A\right)
\right]  _{\infty}=\left[  \overline{\operatorname*{co}}\left(  \cup
_{t\in\mathcal{P}}\partial_{\varepsilon}\left(  \beta_{t}f_{t}\right)  (x)\cup
A\right)  \right]  _{\infty}. \label{tw}%
\end{equation}

\end{lem}

\begin{dem}
To prove (\ref{one}), we first note that $A_{\delta}$ is non-decreasing with
respect to $\delta>0.$ So, we may assume, without loss of generality,\ that
$\delta\geq1$. For each\ $t\in\mathcal{P}$ and $0<\alpha_{t}\leq\beta_{t},$
set\ $\lambda_{t}:=\alpha_{t}/(2\delta\beta_{t})$ $(\in~]0,1/2])$ and
$g_{t}:=\beta_{t}f_{t}.$ Since the sets $\partial_{\varepsilon}\left(
2\lambda_{t}\delta g_{t}\right)  (x)$ and $\partial_{\varepsilon}\left(
2(1-\lambda_{t})\delta g_{t}\right)  (x)$ are nonempty, (\ref{SupDomNega})
entails
\begin{align*}
A_{\delta} &  \subset\left[  \overline{\operatorname*{co}}\left(
\cup_{0<\lambda_{t}\leq1}\left(  \partial_{\varepsilon}\left(  2\lambda
_{t}\delta g_{t}\right)  (x)\cup\partial_{\varepsilon}\left(  2(1-\lambda
_{t})\delta g_{t}\right)  (x)\right)  \cup A\right)  \right]  _{\infty}\\
&  =\left[  \overline{\operatorname*{co}}\left(  \cup_{0<\lambda_{t}\leq
1}\left(  \partial_{\varepsilon}\left(  2\lambda_{t}\delta g_{t}\right)
(x)+\partial_{\varepsilon}\left(  2(1-\lambda_{t})\delta g_{t}\right)
(x)\right)  \cup A\right)  \right]  _{\infty}\\
&  \subset\left[  \overline{\operatorname*{co}}\left(  \cup_{t\in\mathcal{P}%
}\partial_{2\varepsilon}\left(  2\delta g_{t}\right)  (x)\cup A\right)
\right]  _{\infty}\\
&  =\left[  \overline{\operatorname*{co}}\left(  \cup_{t\in\mathcal{P}}%
2\delta\partial_{\varepsilon/\delta}g_{t}(x)\cup A\right)  \right]  _{\infty
}\subset\left[  \overline{\operatorname*{co}}\left(  \cup_{t\in\mathcal{P}%
}\partial_{\varepsilon}g_{t}(x)\cup A\right)  \right]  _{\infty},
\end{align*}
thereby proving (\ref{one}).

To establish (\ref{tw}), define\ $M_{1}:=\sup_{t\in\mathcal{P}}\alpha_{t}%
\beta_{t}^{-1}$ and $M_{2}:=\sup_{t\in\mathcal{P}}\alpha_{t}^{-1}\beta_{t},$
so that $M_{1},M_{2}>0$ and $\alpha_{t}\leq M_{1}\beta_{t}\leq M_{1}%
M_{2}\alpha_{t}$ for all $t\in\mathcal{P}$. Therefore (\ref{tw}) follows
directly from\ (\ref{one}).
\end{dem}

\begin{lem}
\label{lemcone} Assume that $f_{t}\in\Gamma_{0}(X)$ for all $t\in T,$ and
denote $f:=\sup_{t\in T}f_{t}.$ Fix\ $\varepsilon>0,$ $x\in\operatorname*{dom}%
f,$ and $L\in\mathcal{F}(x)$ such that\ $\operatorname*{dom}f_{t}\subset L,$
for all $t\in T.$ Let\ $g_{t}$ denote the restriction of $f_{t}$ to $L.$ Then
$A:={\cup}_{t\in T}\partial_{\varepsilon}f_{t}(x)\neq\emptyset$, $B:={\cup
}_{t\in T}\partial_{\varepsilon}g_{t}(x)\neq\emptyset,$ and we have
\begin{equation}
\left[  \overline{\operatorname*{co}}({B})\right]  _{\infty}\subset\left\{
x_{\mid L}^{\ast}:x^{\ast}\in\left[  \overline{\operatorname*{co}}%
({A})\right]  _{\infty}\right\}  . \label{o1}%
\end{equation}

\end{lem}

\begin{dem}
Take\ $y^{\ast}\in\left[  \overline{\operatorname*{co}}(B)\right]  _{\infty}$
and $x_{0}^{\ast}\in A$, and denote the restriction of $x_{0}^{\ast}$ to
$L$\ by $y_{0}^{\ast}:=x_{0\mid L}^{\ast}$ $(\in B,$ by \cite[Exercise
55(i)]{CHLBook}). Choose a $\theta$-neighborhood $U\subset(X^{\ast},w^{\ast}%
)$, so that $\tilde{U}:=\{u_{\mid L}^{\ast}:u^{\ast}\in U\}$ is a $\theta
$-neighborhood in $L^{\ast}$ (see \cite[page 31]{CHLBook}). Next, by the
extension theorem, there exists\ $x^{\ast}\in X^{\ast}$ such that $y^{\ast
}=x_{\mid L}^{\ast}$ and, so, $y_{0}^{\ast}+\gamma y^{\ast}\in\overline
{\operatorname*{co}}B\subset\operatorname*{co}B+\tilde{U}$ for any\ $\gamma
>0.$ Again \cite[Exercise 55(i)]{CHLBook} yields some\ $z^{\ast}%
\in\operatorname*{co}A$ and $u^{\ast}\in U$ such that $y_{0}^{\ast}+\gamma
y^{\ast}=z_{\mid L}^{\ast}+u_{\mid L}^{\ast}\in z_{\mid L}^{\ast}+\tilde{U}.$
Hence,\
\[
x_{0}^{\ast}+\gamma x^{\ast}\in z^{\ast}+u^{\ast}+L^{\perp}\subset
\operatorname*{co}\left(  {A}+L^{\perp}\right)  +U.
\]
Moreover, since $\operatorname*{dom}f_{t}\subset L$, we have\ $\partial
_{\varepsilon}f_{t}(x)+L^{\perp}=\partial_{\varepsilon}f_{t}(x)+\partial
\mathrm{I}_{L}(x)\subset\partial_{\varepsilon}(f_{t}+\mathrm{I}_{L}%
)(x)=\partial_{\varepsilon}f_{t}(x)$. So, for all $\gamma>0,$ $x_{0}^{\ast
}+\gamma x^{\ast}\in\operatorname*{co}A+U,$ and the arbitrariness of $U$
implies $x_{0}^{\ast}+\gamma x^{\ast}\in\overline{\operatorname*{co}}A.$ Thus,
$x^{\ast}\in\left[  \overline{\operatorname*{co}}A\right]  _{\infty}.$
\end{dem}

\section{Characterization of the normal cone to the domain\label{NormalCone}}

The aim of this\ section is\ to characterize the normal cone to the
(effective) domain of the supremum function
\[
f:=\sup_{t\in T}f_{t},
\]
where each\ $f_{t}:X\rightarrow\overline{\mathbb{R}},$ $t\in T,$ is assumed to
be convex.\ The formula for\ $\mathrm{N}_{\operatorname*{dom}f}(x)$ will rely
on the $\varepsilon$-subdifferential of the functions that are $\varepsilon
$-active at $x,$ as well as on that of the non-$\varepsilon$-active functions,
the latter entering with suitable weights. Nonproper functions are also
relevant, since they contribute through the $\varepsilon$-normal sets
$\mathrm{N}_{\operatorname*{dom}f_{t}}^{\varepsilon}(x).$

\medskip

We begin by introducing a family of weight parameters that penalize those
constraints which are not $\varepsilon$-active at the given nominal point
$x\in\operatorname*{dom}f.$

\begin{definiti}
Let $f_{t}:X\rightarrow\mathbb{R}_{\infty},$ $t\in T,$ be convex, and
denote\ $f:=\sup_{t\in T}f_{t}$.\ For $x\in\operatorname*{dom}f$ and
$\varepsilon>0,$ we define the penalization weights of non-$\varepsilon
$-active indices at $x$ by
\begin{equation}
\rho_{t,\varepsilon}:=\frac{-\varepsilon}{2f_{t}(x)-2f(x)+\varepsilon}\text{
if }t\in T\setminus T_{\varepsilon}(x),\text{ }\rho_{t,\varepsilon}:=1\text{ if\ }t\in
T_{\varepsilon}(x),\label{vetb}%
\end{equation}
where $T_{\varepsilon}(x)$ is the set of $\varepsilon$-active constraints at
$x$ given by
\begin{equation}
T_{\varepsilon}(x):=\{t\in T:f_{t}(x)\geq f(x)-\varepsilon\}.\label{teps2}%
\end{equation}
In case of possible\ confusion, we sometimes denotes $\rho_{t,\varepsilon}$
by\ $\rho_{t,\varepsilon,x}$ to emphasize its dependence on the reference
point $x.$
\end{definiti}

\medskip The next lemma clarifies the role of the weights\ $\rho
_{t,\varepsilon,x}$ introduced above,$\ $showing that they ensure the
quantities\ $\rho_{t,\varepsilon}f_{t}(x)$ are bounded below uniformly in
$t\in T.$

\begin{lem}
\label{lema1} Let $f_{t}:X\rightarrow\mathbb{R}_{\infty},$ $t\in T,$ be
convex, and let $f:=\sup_{t\in T}f_{t}$.\ Fix\ $x\in\operatorname*{dom}f$ and
$\varepsilon>0.$ Then each $\rho_{t,\varepsilon}$ is well-defined and, for all
$t\in T,$
\begin{equation}
0<\rho_{t,\varepsilon}=\varepsilon(\max\{2f(x)-2f_{t}(x)-\varepsilon
,\varepsilon\})^{-1}\leq1, \label{1.1}%
\end{equation}%
\begin{equation}
\min\{f(x),0\}\leq\rho_{t,\varepsilon}f(x)\leq\rho_{t,\varepsilon}%
f_{t}(x)+\varepsilon. \label{4.5}%
\end{equation}
In particular, for all $\delta>0$ and $t\in T_{\varepsilon+\delta}(x),$
\begin{equation}
\varepsilon/(\varepsilon+2\delta)\leq\rho_{t,\varepsilon}\leq1. \label{4.4}%
\end{equation}

\end{lem}

\begin{dem}
We\ denote $\rho_{t}:=\rho_{t,\varepsilon}.$ By replacing each function
$f_{t}$ with $f_{t}-f(x),$ we may\ suppose, without loss of generality, that
$f(x)=0.$ Relation (\ref{1.1}) follows from the fact that\ $t\in T\setminus
T_{\varepsilon}(x)$ is equivalent to\ $-2f_{t}(x)-\varepsilon>-f_{t}%
(x)>\varepsilon$. Moreover, given any $\delta>0$ and\ $t\in T_{\varepsilon
+\delta}(x)\setminus T_{\varepsilon}(x),$ we have $\frac{\varepsilon}{\rho
_{t}}=-2f_{t}(x)-\varepsilon\leq2\delta+\varepsilon,$ so that\ $\rho_{t}%
\geq\varepsilon/(2\delta+\varepsilon).$ Thus, (\ref{4.4}) holds as\ $\rho
_{t}=1\geq\varepsilon/(2\delta+\varepsilon),$ for all\ $t\in T_{\varepsilon
}(x).$ Finally, if\ $t\in T\setminus T_{\varepsilon}(x)$, then
\[
\rho_{t}f_{t}(x)=-\varepsilon\frac{-f_{t}(x)}{-2f_{t}(x)-\varepsilon
}=-\varepsilon\frac{-f_{t}(x)}{(-f_{t}(x)-\varepsilon)-f_{t}(x)}%
\geq-\varepsilon\frac{-f_{t}(x)}{-f_{t}(x)}=-\varepsilon,
\]
and\ (\ref{4.5}) holds. Thus, we are done, since (\ref{4.5}) holds trivially
for\ $t\in T_{\varepsilon}(x)$ (where $\rho_{t}=1).$
\end{dem}

\medskip

The next lemma describes certain property of the weights $\rho_{t,\varepsilon
}$ that will be used later.

\begin{lem}
\label{lema2} Let $f_{t}:X\rightarrow\mathbb{R}_{\infty},$ $t\in T,$ be convex
functions, and set $f:=\sup_{t\in T}f_{t}.$\ Fix\ $x\in f^{-1}(\mathbb{R}).$

$(i)$ Let\ $\varepsilon>0$\ such that $T\setminus T_{\varepsilon}%
(x)\neq\emptyset,$ and define\ $f_{0}:=\sup_{t\in T\setminus T_{\varepsilon
(x)}}f_{t}$ and $\tilde{\rho}_{t,\varepsilon}:=\varepsilon(\max\{2f_{0}%
(x)-2f_{t}(x)-\varepsilon,\varepsilon\})^{-1},$ $t\in T.$ Then
\[
1\leq\inf_{t\in T}\frac{\tilde{\rho}_{t,\varepsilon}}{\rho_{t,\varepsilon}%
}\leq\sup_{t\in T}\frac{\tilde{\rho}_{t,\varepsilon}}{\rho_{t,\varepsilon}%
}<+\infty.
\]

$(ii)$ For $i=1,2,$ let $\varepsilon_{i}>0$\ such that $T\setminus
T_{\varepsilon_{i}(x)}\neq\emptyset,$ and define\ $f_{i}:=\sup_{t\in
T\setminus T_{\varepsilon_{i}(x)}}f_{t}$ and $\tilde{\rho}_{t,\varepsilon_{i}%
}:=\varepsilon(\max\{2f_{i}(x)-2f_{t}(x)-\varepsilon_{i},\varepsilon
_{i}\})^{-1},$ $t\in T.$ Then
\[
\sup_{t\in T}\frac{\tilde{\rho}_{t,\varepsilon_{1}}}{\tilde{\rho
}_{t,\varepsilon_{2}}}<+\infty,
\]
and the family of sets $B_{\varepsilon}:=\overline{\operatorname*{co}}%
(\cup_{t\in\mathcal{P}\setminus T_{\varepsilon}(x)}\partial_{\varepsilon
}(\tilde{\rho}_{t,\varepsilon}f_{t})(x)),$ $\varepsilon>0,$ is nondecreasing
as $\varepsilon\downarrow0.$
\end{lem}

\begin{dem}
$(i)$ Since $f\geq f_{0}$ and
\begin{equation}
\frac{\tilde{\rho}_{t,\varepsilon}}{\rho_{t,\varepsilon}}=\frac{\max
\{2f(x)-2f_{t}(x)-\varepsilon,\varepsilon\}}{\max\{2f_{0}(x)-2f_{t}%
(x)-\varepsilon,\varepsilon\}},\text{ for all }t\in T, \label{df}%
\end{equation}
it follows that $\inf_{t\in T}\frac{\tilde{\rho}_{t,\varepsilon}}%
{\rho_{t,\varepsilon}}\geq1.$ Next, choose\ $M_{0}>2(\varepsilon
-f(x))\ $sufficiently large so that, for every real $M$ with $\left\vert
M\right\vert \geq M_{0},$ we have $\frac{\max\{2f(x)-M-\varepsilon
,\varepsilon\}}{\max\{2f_{0}(x)-M-\varepsilon,\varepsilon\}}\leq2.$ If\ $t\in
T$ is such that $-2f_{t}(x)\geq M_{0},$ then (\ref{df}) yields $\frac
{\tilde{\rho}_{t,\varepsilon}}{\rho_{t,\varepsilon}}\leq2.$ Otherwise,
$-2f_{t}(\theta)<M_{0}$ and we have $\frac{\tilde{\rho}_{t,\varepsilon}}%
{\rho_{t,\varepsilon}}\leq\frac{\max\{2f(x)+M_{0}-\varepsilon,\varepsilon
\}}{\max\{2f_{0}(x)-2f_{t}(x)-\varepsilon,\varepsilon\}}\leq\frac
{2f(x)+M_{0}-\varepsilon}{\varepsilon},$ thereby proving that $\sup_{t\in
T}\frac{\tilde{\rho}_{t,\varepsilon}}{\rho_{t,\varepsilon}}<+\infty.$

$(ii)$ Choose\ $M_{1}>2(\varepsilon_{2}-\max_{i=1,2}f_{i}(x))\ $sufficiently
large so that, for every real $M$ with $\left\vert M\right\vert \geq M_{1},$
we have $\frac{\max\{2f_{2}(x)-M-\varepsilon_{2},\varepsilon_{2}\}}%
{\max\{2f_{1}(x)-M-\varepsilon_{1},\varepsilon_{1}\}}\leq2+\varepsilon
_{2}/\varepsilon_{1}.$ If\ $t\in T$ is such that $-2f_{t}(x)\geq M_{1},$ then
\[
\frac{\tilde{\rho}_{t,\varepsilon_{1}}}{\tilde{\rho}_{t,\varepsilon_{2}}%
}=\frac{\varepsilon_{1}\max\{2f_{2}(x)-2f_{t}(x)-\varepsilon_{2}%
,\varepsilon_{2}\}}{\varepsilon_{2}\max\{2f_{2}(x)-2f_{t}(x)-\varepsilon
_{1},\varepsilon_{1}\}}\leq2+\varepsilon_{2}/\varepsilon_{1}.
\]
Otherwise, $-2f_{t}(\theta)<M_{1}$ and we have $$\frac{\tilde{\rho
}_{t,\varepsilon_{1}}}{\tilde{\rho}_{t,\varepsilon_{2}}}\leq\frac{\max
\{2f_{2}(x)+M_{1}-\varepsilon_{2},\varepsilon_{2}\}}{\max\{2f_{1}%
(x)-2f_{t}(x)-\varepsilon_{1},\varepsilon_{1}\}}\leq\frac{2f_{2}%
(x)+M_{1}-\varepsilon_{2}}{\varepsilon_{1}},$$ thereby proving that $\sup_{t\in
T}\frac{\tilde{\rho}_{t,\varepsilon_{1}}}{\tilde{\rho}_{t,\varepsilon_{2}}%
}<+\infty.$

Assume now that $\varepsilon_{1}\leq\varepsilon_{2}.$ By the previous
assertion, there exist $r_{1},$ $r_{2}>0$ such that $\varepsilon_{1}%
\tilde{\rho}_{t,\varepsilon_{2}}\leq r_{1}\varepsilon_{2}\tilde{\rho
}_{t,\varepsilon_{1}}\leq r_{1}r_{2}\varepsilon_{1}\tilde{\rho}_{t,\varepsilon
_{2}},$ for all $t\in T.$ Hence, applying Lemma \ref{ChangeParameter} yields\
\[
\left[  B_{\varepsilon_{2}}\right]  _{\infty}=\left[  \underset{t\in
\mathcal{P}\setminus T_{\varepsilon_{2}}(x)}{\cup}\partial_{\varepsilon_{2}%
}(\tilde{\rho}_{t,\varepsilon_{1}}f_{t})(x)\right]  _{\infty}\subset\left[
B_{\varepsilon_{1}}\right]  _{\infty},
\]
where in the last inclusion we used the easy fact that $\mathcal{P}\setminus
T_{\varepsilon_{2}}(x)\subset\mathcal{P}\setminus T_{\varepsilon_{1}}(x).$
\end{dem}

\medskip The following lemma is\ needed in\ the proof of Theorem
\ref{Theo13b}. Its purpose is to\ control the number of constraints that enter
into the subdifferential of the supremum, relative to the dimension of the
underlying (finite-dimensional) subspace.

\begin{lem}
\label{SubFiniteCaseb} Given $f_{1},\cdots,f_{m}\in\Gamma_{0}(\mathbb{R}^{n})$
with $m>n+1$, and\ $f:=\max_{1\leq i\leq m}f_{i},$ for every $x\in
\operatorname*{dom}f$ and $\varepsilon>0$ we have
\begin{equation}
\partial_{\varepsilon}f(x)\subset\left\{  \partial_{\varepsilon}f_{\lambda
}(x):\lambda\in\Delta_{m},\text{ }f_{\lambda}(x)\geq f(x)-\varepsilon
,\ \left\vert \operatorname{supp}\lambda\right\vert \leq n+1\right\}  ,
\label{ad}%
\end{equation}
where $\left\vert \operatorname{supp}\lambda\right\vert $ denotes the cardinal
of $\operatorname{supp}\lambda$ and $f_{\lambda}:=\sum_{i\in
\operatorname{supp}\lambda}\lambda_{i}f_{i}$.
\end{lem}

\begin{dem}
First, assume that $\theta\in\partial_{\varepsilon}f(x).$ According
to\ \cite[Theorem 1]{Dr82} (see, also, \cite[Corollary 5]{CH24a}), there
exists\ $\lambda\in\Delta_{m}$ with\ $\left\vert \operatorname{supp}%
\lambda\right\vert \leq n+1$ such that
\[
f(x)\leq\inf_{\mathbb{R}^{n}}f+\varepsilon=\inf_{\mathbb{R}^{n}}f_{\lambda
}+\varepsilon.
\]
Thus, $f_{\lambda}(x)\leq f(x)\leq\inf_{\mathbb{R}^{n}}f_{\lambda}%
+\varepsilon\leq f_{\lambda}(x)+\varepsilon,$ and we deduce that $\theta
\in\partial_{\varepsilon}f_{\lambda}(x)$ with $f_{\lambda}(x)\geq
f(x)-\varepsilon.$ More generally, given\ $x^{\ast}\in\partial_{\varepsilon
}f(x),$ we consider the functions\ $g_{i}:=f_{i}-\left\langle x^{\ast}%
,\cdot\right\rangle \in\Gamma_{0}(\mathbb{R}^{n}),$ $i=1,\cdots,m,$ and
$g:=\max_{1\leq i\leq m}g_{i}.$ Then $\theta\in\partial_{\varepsilon}g(x),$
and the arguments above yield a $\lambda\in\Delta_{m}$ such that $\left\vert
\operatorname{supp}\lambda\right\vert \leq n+1,$ $\sum_{i\in
\operatorname{supp}\lambda}\lambda_{i}g_{i}(x)\geq g(x)-\varepsilon$ and
$\theta\in\partial_{\varepsilon}(\sum_{i\in\operatorname{supp}\lambda}%
\lambda_{i}g_{i})(x)=\partial_{\varepsilon}f_{\lambda}(x)-x^{\ast}.$
Equivalently, $\lambda$ satisfies $f_{\lambda}(x)\geq f(x)-\varepsilon$ and
$x^{\ast}\in\partial_{\varepsilon}f_{\lambda}(x),$ showing that\ (\ref{ad}) holds.
\end{dem}

\medskip The following\ result provides a characterization of\ the normal
cone\ $\mathrm{N}_{\operatorname*{dom}f}(x)$ in the case where all
functions\ $f_{t}$ are in\ $\Gamma_{0}(X)$. More general versions of this
result are obtained later as corollaries. Remember that the weights\ $\rho
_{t,\varepsilon,x}$ are\ defined in (\ref{vetb}).

\begin{theo}
\label{Theo13b} Let $f_{t}:X\rightarrow\mathbb{R}_{\infty},\ t\in T,$
be\ proper lsc convex\ functions, and denote\ $f:=\sup_{t\in T}f_{t}.$ Then,
for every $x\in\operatorname*{dom}f,$ we have
\begin{equation}
\mathrm{N}_{\operatorname*{dom}f}(x)=\left[  \overline{\operatorname*{co}%
}\left(  {%
{\textstyle\bigcup\limits_{t\in T}}
}\partial_{\varepsilon}(\alpha_{t}f_{t})(x)\right)  \right]  _{\infty},
\label{mainfb}%
\end{equation}
for any\ $\varepsilon>0$ and any $\alpha_{t}\geq\rho_{t,\varepsilon}$ such
that $\sup_{t\in T}\alpha_{t}<+\infty$ and $\inf_{t\in T}\alpha_{t}%
f_{t}(x)>-\infty.$ Consequently,
\begin{equation}
\mathrm{N}_{\operatorname*{dom}f}(x)=%
{\textstyle\bigcap\limits_{\varepsilon>0}}
\overline{\operatorname*{co}}\left(  {%
{\textstyle\bigcup\limits_{t\in T}}
}\partial_{\varepsilon}(\varepsilon\alpha_{t}f_{t})(x)\right)  ,
\label{mainfc}%
\end{equation}
provided that the latter set is nonempty; otherwise, $\mathrm{N}%
_{\operatorname*{dom}f}(x)=\{\theta\}.$
\end{theo}

\begin{rem}
[before the proof]\label{rr}More precisely, as shown in the proof, the
inclusion \textquotedblleft$\subset$\textquotedblright\ in (\ref{mainfb}%
)\ holds without assuming the condition $\inf_{t\in T}\alpha_{t}%
f_{t}(x)>-\infty;$ in fact, this last assumption is required only for the
opposite inclusion.
\end{rem}

\begin{dem}
Fix\ $x\in\operatorname*{dom}f$ and $\varepsilon>0$. Since each\ $f_{t}$ is
proper, we have $f(x)\in\mathbb{R}$, and we may assume without loss of
generality\ that $f(x)=0.$ We divide the proof into five\ steps. In Steps 1-3,
we prove the inclusion\ \textquotedblleft$\subset$\textquotedblright\ in
(\ref{mainfb}) under the assumption that the parameters\ $\alpha_{t}$
satisfy\ $\alpha_{t}\geq\rho_{t,\varepsilon}$ and\ $\sup_{t\in T}\alpha
_{t}<+\infty.$ The additional condition $\inf_{t\in T}\alpha_{t}%
f_{t}(x)>-\infty$ is imposed only in Step 4, where we prove the inclusion
\textquotedblleft$\supset$\textquotedblright\ in (\ref{mainfb}). Finally, the
proof of (\ref{mainfc}) is given in Step 5.

\textbf{Step 1.} We suppose in this\ first step that $X$ is (finite)\ $n$%
-dimensional. According to \cite[Theorem 7]{HaLo22}, for any fixed index
$t_{0}\in T_{\varepsilon/2}(x)$ we have
\begin{equation}
\mathrm{N}_{\operatorname*{dom}f}(x)=\left[  \overline{\operatorname*{co}%
}\left(  {\cup}_{J\in\mathcal{T}}\partial_{\varepsilon/2}f_{J}(x)\right)
\right]  _{\infty}, \label{tt}%
\end{equation}
where $\mathcal{T}:=\{J\subset T:t_{0}\in J,$ $\left\vert J\right\vert
<+\infty\}$ and $f_{J}:=\max_{t\in J}f_{t},$ $J\in\mathcal{T}.$ More
precisely, for each\ fixed\ $J:=\{t_{0},$ $t_{1},$ $\cdots,$ $t_{m-1}%
\}\in\mathcal{T}$ with $m\geq n+1,$\ we have $f_{J}(x)\geq f_{t_{0}}%
(x)\geq-\varepsilon/2.$ Moreover, if we set $f_{\lambda}:=\sum_{t\in
\operatorname{supp}\lambda}\lambda_{t}f_{t},$ then\ applying Lemma
\ref{SubFiniteCaseb} gives
\begin{equation}
\partial_{\frac{\varepsilon}{2}}f_{J}(x)\subset{%
{\textstyle\bigcup_{\substack{\lambda\in\Delta_{m}\\\operatorname{supp}%
\lambda\subset J}}}
}\left\{  \partial_{\frac{\varepsilon}{2}}f_{\lambda}(x):f_{\lambda}(x)\geq
f_{J}(x)-\frac{\varepsilon}{2}\geq-\varepsilon,\text{ }\left\vert
\operatorname{supp}\lambda\right\vert \leq n+1\right\}  . \label{eq1}%
\end{equation}
We now claim that any $\lambda\in\Delta_{m}$ appearing in (\ref{eq1})
satisfies
\begin{equation}
\lambda_{t}\leq2\alpha_{t},\text{ for all }t\in\operatorname{supp}\lambda,
\label{claim}%
\end{equation}
where $\alpha_{t}\geq\rho_{t,\varepsilon}$ are such that\ $\sup_{t\in T}%
\alpha_{t}<+\infty.$ To prove this claim we consider the set $J_{\varepsilon
}:=\{t\in\operatorname{supp}\lambda:f_{t}(x)\geq-\varepsilon\},$ which is
nonempty; otherwise, $f_{\lambda}(x)<-\varepsilon,$ and we get a contradiction
with the choice of $\lambda$ in (\ref{eq1}). We observe that
\begin{equation}%
{\textstyle\sum\nolimits_{t\in J_{\varepsilon}}}
\lambda_{t}(f_{t}(x)+\varepsilon/2)+%
{\textstyle\sum\nolimits_{t\in(\operatorname{supp}\lambda)\setminus
J_{\varepsilon}}}
\lambda_{t}(f_{t}(x)+\varepsilon/2)=f_{\lambda}(x)+\varepsilon/2\geq
-\varepsilon/2, \label{IneSumft}%
\end{equation}
and, since\ $-\varepsilon\leq f_{t}(x)\leq f(x)=0$ for every $t\in
J_{\varepsilon},$
\begin{equation}
\varepsilon/2\geq(\varepsilon/2)%
{\textstyle\sum\nolimits_{t\in J_{\varepsilon}}}
\lambda_{t}\geq%
{\textstyle\sum\nolimits_{t\in J_{\varepsilon}}}
\lambda_{t}(f_{t}(x)+\varepsilon/2)\geq-(\varepsilon/2)%
{\textstyle\sum\nolimits_{t\in J_{\varepsilon}}}
\lambda_{t}. \label{com}%
\end{equation}
Consequently, combining\ (\ref{IneSumft}) and\ (\ref{com}) entails\
\[%
{\textstyle\sum\nolimits_{t\in(\operatorname{supp}\lambda)\setminus
J_{\varepsilon}}}
\lambda_{t}(f_{t}(x)+\varepsilon/2)\geq-(\varepsilon/2)-%
{\textstyle\sum\nolimits_{t\in J_{\varepsilon}}}
\lambda_{t}(f_{t}(x)+\varepsilon)\geq-\varepsilon.
\]
At the same time, each\ index $t\in(\operatorname{supp}\lambda)\setminus
J_{\varepsilon}$ belongs to $T\setminus T_{\varepsilon}(x)$ and
satisfies\ $\lambda_{t}(f_{t}(x)+\varepsilon/2)<0.$ Then (\ref{com}) further
implies\ $$\lambda_{t}(f_{t}(x)+\varepsilon/2)\geq%
{\textstyle\sum\nolimits_{s\in(\operatorname{supp}\lambda)\setminus
J_{\varepsilon}}}
\lambda_{s}(f_{s}(x)+\varepsilon/2)\geq-\varepsilon$$ and, consequently,
\[
\lambda_{t}\leq\frac{-\varepsilon}{f_{t}(x)+\varepsilon/2}=\frac
{-2\varepsilon}{2f_{t}(x)+\varepsilon}=2\rho_{t,\varepsilon}\leq2\alpha_{t}.
\]
Therefore the claim (\ref{claim}) is established\ as\ $\lambda_{t}\leq
1\leq2=2\rho_{t,\varepsilon}=2\alpha_{t},$ for all $t\in J_{\varepsilon}$
$(\subset T_{\varepsilon}(x)).$

\textbf{Step 2. }We now apply the Hiriart-Urruty \& Phelps Theorem (see
\cite[Proposition 4.1.16]{CHLBook}) to the functions\ $f_{t}\in\Gamma_{0}(X),$
defined on $X$ (supposed of dimension $n$). For each $J\in\mathcal{T}$,
combining (\ref{eq1}) and (\ref{claim}) entails\
\begin{align*}
  \partial_{\frac{\varepsilon}{2}}f_{J}(x)
  &\subset{%
{\textstyle\bigcup_{\substack{\lambda\in\Delta_{m}\\\operatorname{supp}%
\lambda\subset J}}}
}\left\{
{\textstyle\sum\nolimits_{t\in\operatorname{supp}\lambda}}
\partial_{\varepsilon/2}\left(  \lambda_{t}f_{t}\right)  (x):\left\vert
\operatorname{supp}\lambda\right\vert \leq n+1,\text{ }\lambda_{t}\leq
2\alpha_{t}\right\}  +\varepsilon\mathbb{B}_{X^{\ast}}\\
&  \subset \varepsilon\mathbb{B}%
_{X^{\ast}}+ (n+1)\operatorname*{co}\{\cup_{t\in\operatorname{supp}\lambda,\text{
}\lambda\in\Delta_{m},\lambda\subset J}\{\partial_{\varepsilon/2}\left(
\lambda_{t}f_{t}\right)  (x):
\\
& \qquad   \qquad  \qquad  \qquad  \qquad  \qquad  \qquad  \qquad \qquad  \qquad \left\vert \operatorname{supp}\lambda\right\vert
\leq n+1,\text{ }\lambda_{t}\leq2\alpha_{t}\}\},
\end{align*}
where $\mathbb{B}_{X^{\ast}}$ denotes the unit ball in the ($n$-dimensional)
dual space $X^{\ast}$ of $X$. Therefore, using \cite[(3.52), p. 80]{CHLBook},
(\ref{tt}) simplifies to\
\begin{align*}
\mathrm{N}_{\operatorname*{dom}f}(x)  &  =\left[  \overline{\operatorname*{co}%
}\left(  {\cup}_{J\in\mathcal{T}}\partial_{\varepsilon/2}f_{J}(x)\right)
\right]  _{\infty}=\left[  \overline{\operatorname*{co}}\left(  {\cup}%
_{J\in\mathcal{T}}(n+1)^{-1}\partial_{\varepsilon/2}f_{J}(x)\right)  \right]
_{\infty}\\
&  \subset\left[  \overline{\operatorname*{co}}\left(  \cup_{0<\lambda_{t}%
\leq2\alpha_{t}}\partial_{\varepsilon}\left(  \lambda_{t}f_{t}\right)
(x)+\varepsilon\mathbb{B}_{X^{\ast}}\right)  \right]  _{\infty}=\left[
\overline{\operatorname*{co}}\left(  \cup_{0<\lambda_{t}\leq2\alpha_{t}%
}\partial_{\varepsilon}\left(  \lambda_{t}f_{t}\right)  (x)\right)  \right]
_{\infty},
\end{align*}
and Lemma \ref{ChangeParameter} entails the inclusion \textquotedblleft%
$\subset$\textquotedblright\ in (\ref{mainfb}).

\textbf{Step 3. }To show that the inclusion \textquotedblleft$\subset
$\textquotedblright\ in (\ref{mainfb}) also holds when $X$ is
infinite-dimensional, we choose a (weak*-)\ $\theta$-neighborhood $U\subset
X^{\ast},$\ and select $L\in\mathcal{F}(x)$ such that $L^{\perp}\subset U$.
Let\ $g_{t}\in\Gamma_{0}(L),$ $t\in T,$ be the restriction of $f_{t}%
+\mathrm{I}_{L}$ to the finite-dimensional space\ $L,$ and denote
$g:=\sup_{t\in T}g_{t}=f+\mathrm{I}_{L}$. Since $x\in L,$ the families
$\{f_{t},$ $t\in T\}$ and $\{g_{t},$ $t\in T\}$ define the same $\varepsilon
$-active set $T_{\varepsilon}(x)$ and the same numbers $\rho_{t,\varepsilon}$
(see (\ref{vetb})). Then, since the inclusion \textquotedblleft$\subset
$\textquotedblright\ in (\ref{mainfb}) holds for\ the family $\{g_{t},$ $t\in
T\},$ we obtain $\mathrm{N}_{\operatorname*{dom}g}(x)\subset\left[
\overline{\operatorname*{co}}\left(  \cup_{t\in T}\partial_{\varepsilon
}(\alpha_{t}g_{t})(x)\right)  \right]  _{\infty}.$ Moreover, knowing that
$\mathrm{N}_{\operatorname*{dom}g}(x)=\{x_{\mid L}^{\ast}\in L^{\ast}:x^{\ast
}\in\mathrm{N}_{L\cap\operatorname*{dom}f}(x)\}$ (see \cite[Exercise
55(i)]{CHLBook}) and
\[
\left[  \overline{\operatorname*{co}}\left(  \cup_{t\in T}\partial
_{\varepsilon}(\alpha_{t}g_{t})(x)\right)  \right]  _{\infty}\subset\left\{
x_{\mid L}^{\ast}:x^{\ast}\in\left[  \overline{\operatorname*{co}}\left(
\cup_{t\in T}\partial_{\varepsilon}(\alpha_{t}f_{t}+\mathrm{I}_{L})(x)\right)
\right]  _{\infty}\right\}  ,
\]
by Lemma \ref{lemcone}, we conclude that
\begin{align*}
\mathrm{N}_{L\cap\operatorname*{dom}f}(x)  &  \subset\left[  \overline
{\operatorname*{co}}\left(  \cup_{t\in T}\partial_{\varepsilon}(\alpha
_{t}f_{t}+\mathrm{I}_{L})(x)\right)  \right]  _{\infty}+L^{\perp}\\
&  \subset\left[  \overline{\operatorname*{co}}\left(  \cup_{t\in T}%
\partial_{\varepsilon}(\alpha_{t}f_{t})(x)+L^{\perp}\right)  \right]
_{\infty}+L^{\perp}
\\
&\subset
\left[  \overline{\operatorname*{co}}\left(
\cup_{t\in T}\partial_{\varepsilon}(\alpha_{t}f_{t})(x)+L^{\perp}\right)
\right]  _{\infty},
\end{align*}
where the second and third inclusions come from \cite[Proposition
4.1.16]{CHLBook} and (\ref{su}), respectively. So, by intersecting over all
finite-dimensional subspaces $L^{\prime}\supset L$, we obtain
\[
\mathrm{N}_{\operatorname*{dom}f}(x)\subset\cap_{L^{\prime}\supset
L}\mathrm{N}_{L^{\prime}\cap\operatorname*{dom}f}(x)\subset\cap_{L^{\prime
}\supset L}\left[  \overline{\operatorname*{co}}\left(  \cup_{t\in T}%
\partial_{\varepsilon}(\alpha_{t}f_{t})(x)+L^{\perp}\right)  \right]
_{\infty}.
\]
Finally, the inclusion \textquotedblleft$\subset$\textquotedblright\ in
(\ref{mainfb}) follows from\ (\ref{lemIntReca}).

\textbf{Step 4. }To prove the inclusion \textquotedblleft$\supset
$\textquotedblright\ in (\ref{mainfb}), fix\ $\varepsilon>0$ and let $\alpha$
be as in the statement of the theorem. Choose $t_{0}\in T$ and take\ $x_{0}%
^{\ast}\in\partial_{\varepsilon}(\alpha_{t_{0}}f_{t_{0}})(x)$ (which\ is
nonempty because\ $f_{t}\in\Gamma_{0}(X)$ and $\alpha_{t_{0}}>0$). Given
$x^{\ast}$ in the right hand side of (\ref{mainfb}) and $\beta>0,$\ we find
nets $(k_{j})_{j}\subset\mathbb{N}$, $(\lambda_{j,1},\cdots,\lambda_{j,k_{j}%
})_{j}\subset\Delta_{k_{j}},$ $(t_{j,1},\cdots,t_{j,k_{j}})_{j}\subset T$ and
elements $x_{j,i}^{\ast}\in\partial_{\varepsilon}(\alpha_{t_{j,i}}f_{t_{j,i}%
})(x)$ for $i=1,\cdots,k_{j}$ such that\ $x_{0}^{\ast}+\beta x^{\ast}=\lim
_{j}(\lambda_{j,1}x_{j,1}^{\ast}+\cdots+\lambda_{j,k_{j}}x_{j,k_{j}}^{\ast}).$
Fix $y\in\operatorname{dom}f\ (\subset\operatorname{dom}f_{t})$. Since
$\sup_{t\in T}\alpha_{t}$ is finite and $\inf_{t\in T}\alpha_{t}%
f_{t}(x)>-\infty,$ we obtain
\begin{align*}
\langle x_{0}^{\ast}+\beta x^{\ast},y-x\rangle &  =\lim_{j}\langle
\lambda_{j,1}x_{j,1}^{\ast}+\cdots+\lambda_{j,k_{j}}x_{j,k_{j}}^{\ast
},y-x\rangle\\
&  \leq\limsup_{j}%
{\textstyle\sum\nolimits_{i=1}^{k_{j}}}
\lambda_{j,i}\left(  \alpha_{t_{j,i}}f_{t_{j,i}}(y)-\alpha_{t_{j,i}}%
f_{t_{j,i}}(x)+\varepsilon\right) \\
&  \leq\limsup_{j}%
{\textstyle\sum\nolimits_{i=1}^{k_{j}}}
\lambda_{j,i}\left(  \alpha_{t_{j,i}}f^{+}(y)-\alpha_{t_{j,i}}f_{t_{j,i}%
}(x)+\varepsilon\right) \\
&  \leq(\sup_{t\in T}\alpha_{t})f^{+}(y)-\inf_{t\in T}\alpha_{t}%
f_{t}(x)+\varepsilon.
\end{align*}
Thus,$\ $dividing the above inequality by\ $\beta$ and letting $\beta
\rightarrow+\infty$ yields\ $\langle x^{\ast},y-x\rangle\leq0,\ $for all
$y\in\operatorname{dom}f=\operatorname{dom}f^{+},$ which proves that $x^{\ast
}\in\mathrm{N}_{\operatorname{dom}f}(x).$ The proof of the theorem is complete.

\textbf{Step 5.} Set $A_{\varepsilon}:=\overline{\operatorname*{co}}\left(
{\cup}_{t\in T}\partial_{\varepsilon}(\varepsilon\alpha_{t}f_{t})(x)\right)
,$ $\varepsilon>0.$ Since $\lambda A_{\varepsilon}=A_{\lambda\varepsilon}$ for
all $\varepsilon,\lambda>0,$ it follows that $\{\theta\}\cup(\cap
_{\varepsilon>0}A_{\varepsilon})$ is a nonempty closed cone. Thus, taking into
account\ Lemma \ref{ChangeParameter} and (\ref{ricone}), by taking the
intersection over $\varepsilon>0$ in (\ref{mainfb}) we obtain\
\[
\mathrm{N}_{\operatorname*{dom}f}(x)=\cap_{\varepsilon>0}\left[
A_{\varepsilon}\cup\{\theta\}\right]  _{\infty}=\left[  \{\theta\}\cup
\cap_{\varepsilon>0}A_{\varepsilon}\right]  _{\infty}=\{\theta\}\cup
(\cap_{\varepsilon>0}A_{\varepsilon}),
\]
thereby proving (\ref{mainfc}).
\end{dem}

The following corollary provides a typical illustration of Theorem
\ref{Theo13b}, describing the normal cone to an infinite intersection of
closed convex sets without imposing any interiority condition. Representations
involving exact normal cones, rather than the $\varepsilon$-normal sets,
generally require additional interiority conditions and closure-type\ dual
conditions (see, e.g., \cite{DiGoLo06, BJ04}). A finite-dimensional
version\ of this result can be found in \cite[Proposition 3.1]{Ha06}. The
following\ result can also be derived from Corollary \ref{bronds} in the
forthcoming section. 

\begin{cor}
Let $C_{t}\subset X,$ $t\in T,$ be\ closed and convex sets, and denote
$C:=\cap_{t\in T}C_{t}$. Then, for every $x\in C,$
\begin{equation}
\mathrm{N}_{C}(x)=\left[  \overline{\operatorname*{co}}\left(
{\textstyle\bigcup\limits_{t\in T}}
\mathrm{N}_{C_{t}}^{\varepsilon}(x)\right)  \right]  _{\infty},\text{ for
every }\varepsilon>0,\label{ew}%
\end{equation}
and, consequently,
\[
\mathrm{N}_{C}(x)=%
{\textstyle\bigcap\limits_{\varepsilon>0}}
\overline{\operatorname*{co}}\left(
{\textstyle\bigcup\limits_{t\in T}}
\mathrm{N}_{C_{t}}^{\varepsilon}(x)\right)  .
\]
In particular, if\ the $C_{t}$'s are additionally linear subspaces, then
$C^{\perp}=\overline{\operatorname*{span}}\left(  \cup_{t\in T}C_{t}^{\perp
}\right)  $.
\end{cor}

\begin{dem}
The first equality follows by applying Theorem \ref{Theo13b} to the family
$\{f_{t}:=\mathrm{I}_{C_{t}},$ $t\in T\}.$ Next, since $\theta\in
\mathrm{N}_{C_{t}}^{\varepsilon}(x)$ for all $t\in T$ and all $\varepsilon>0,$
taking the intersection over $\varepsilon>0$ in (\ref{ew}) together with
(\ref{ricone}) yields
\[
\mathrm{N}_{C}(x)=\cap_{\varepsilon>0}\left[  \overline{\operatorname*{co}%
}\left(  \cup_{t\in T}\mathrm{N}_{C_{t}}^{\varepsilon}(x)\right)  \right]
_{\infty}=\left[  \cap_{\varepsilon>0}\overline{\operatorname*{co}}\left(
\cup_{t\in T}\mathrm{N}_{C_{t}}^{\varepsilon}(x)\right)  \right]  _{\infty}.
\]
Thus, the second conclusion follows from the fact that\ $\cap_{\varepsilon
>0}\overline{\operatorname*{co}}\left(  \cup_{t\in T}\mathrm{N}_{C_{t}%
}^{\varepsilon}(x)\right)  $ is a closed convex cone, as $\lambda
\mathrm{N}_{C_{t}}^{\varepsilon}(x)=\mathrm{N}_{C_{t}}^{\lambda\varepsilon
}(x),$ for all $\lambda>0$. 
\end{dem}

\smallskip The situation where $\inf_{t\in T}f_{t}(x)>-\infty$, covering the
case where $T$ is finite, leads to a simple\ characterization of the normal
cone. More precisely, we have:

\begin{cor}
\label{corf}Let $f,$ $f_{t}:X\rightarrow\mathbb{R}_{\infty},\ t\in T,$ be as
in Theorem \ref{Theo13b}. Then, for every $x\in\operatorname*{dom}f,$ we have
\[
\mathrm{N}_{\operatorname*{dom}f}(x)\subset\left[  \overline
{\operatorname*{co}}\left(  {%
{\textstyle\bigcup\limits_{t\in T}}
}\partial_{\varepsilon}f_{t}(x)\right)  \right]  _{\infty},\text{ for all
}\varepsilon>0.
\]
Moreover, if\ $\inf_{t\in T}f_{t}(x)>-\infty,$ then (\ref{mainfb}) and
(\ref{mainfc}) hold with $\alpha_{t}=1,$ for all $t\in T.$
\end{cor}

\begin{dem}
Apply Theorem \ref{Theo13b} and Remark \ref{rr} with $\alpha_{t}=1,$ $t\in T.$
\end{dem}

\smallskip The inclusion in the above corollary can still be\ of independent
interest, particularly for deriving necessary optimality conditions. However,
the weights $\alpha_{t}$ are essential\ to obtain an exact representation of
the normal cone, since\ the opposite\ inclusion may fail in general, as shown
by the following one-dimensional example.

\begin{exam}
\label{exam1} Consider the (affine) functions $f_{t}\in\Gamma_{0}(\mathbb{R})$
defined as
\[
f_{t}(x):=tx-t,\text{ }t\in T:=\mathbb{R}_{+},\text{ and }f:=\sup
\nolimits_{t\in T}f_{t}=\mathrm{I}_{]-\infty,1]}.
\]
Then $f(1)=f(0)=0$ and $\operatorname*{dom}f=\left]  -\infty,1\right]  .$
Given\ $\varepsilon>0$ and $x\leq1,$ we have\ $T_{\varepsilon}%
(0)=[0,\varepsilon],$ $\partial_{\varepsilon}f_{t}(x)=\{t\}$ for all\ $t\geq
0$, and so\
\[
\left[  \overline{\operatorname*{co}}\left(  {\cup}_{t\geq0}\partial
_{\varepsilon}f_{t}(0)\right)  \right]  _{\infty}=\left[  \overline
{\operatorname*{co}}\left(  {\cup}_{t\geq0}\partial_{\varepsilon}%
f_{t}(1)\right)  \right]  _{\infty}=\mathbb{R}_{+}.
\]
At the same time, we have
\[
\mathrm{N}_{\operatorname*{dom}f}(1)=\mathbb{R}_{+}=\left[  \overline
{\operatorname*{co}}\left(  {\cup}_{t\geq0}\partial_{\varepsilon}%
f_{t}(1)\right)  \right]  _{\infty},\text{ }\mathrm{N}_{\operatorname*{dom}%
f}(0)=\left\{  0\right\}  \varsubsetneq\mathbb{R}_{+}\left[  \overline
{\operatorname*{co}}\left(  \cup_{t\geq0}\partial_{\varepsilon}f_{t}%
(0)\right)  \right]  _{\infty},
\]
and (\ref{mainfb}) in Theorem \ref{Theo13b} with $\alpha_{t}=1,$ $t\in T,$
is\ applicable at $x=1$ but not at $x=0.$ In fact, we have\
\[
\inf_{t\in T}f_{t}(1)=0(>-\infty)\text{ \ and \ }\inf_{t\in T}f_{t}%
(0)=\inf_{t\geq0}(-t)=-\infty;
\]
that is, the parameters $\alpha_{t}\equiv1$ satisfy the requirements of
Theorem \ref{Theo13b} at $x=1,$ but not at $x=0.$ 

To find the exact representation of $\mathrm{N}_{\operatorname*{dom}f}(0)$ via
Theorem \ref{Theo13b}, we need to consider the parameters $\rho_{t,\varepsilon
,0}$ defined in (\ref{vetb}). Indeed, in the present case we have\ $\rho
_{t,\varepsilon,0}=\varepsilon(\max\{2t-\varepsilon,\varepsilon\})^{-1}$ and
Corollary \ref{corlaylaha} reads
\begin{align}
\mathrm{N}_{\operatorname*{dom}f}(0)&=\left[  \overline{\operatorname*{co}%
}\left(  \underset{t\leq\varepsilon}{\cup}\partial_{\varepsilon}f_{t}(0){\cup
}\underset{t>\varepsilon}{\cup}\partial_{\varepsilon}(\frac{\varepsilon
t}{2t-\varepsilon}f_{t})(0)\right)  \right]  _{\infty}
\\
&=\left[  \overline
{\operatorname*{co}}\left(  [0,\varepsilon]{\cup}\underset{t>\varepsilon}%
{\cup}\left\{  \frac{\varepsilon t}{2t-\varepsilon}\right\}  \right)  \right]
_{\infty}.
\end{align}
Hence,\ since the set $\left\{  \cup_{t>\varepsilon}\varepsilon
t(2t-\varepsilon)^{-1}\right\}  $ is bounded, we deduce that $\mathrm{N}%
_{\operatorname*{dom}f}(0)=\left\{  0\right\}  .$ 
\end{exam}

\medskip Theorem \ref{Theo13b} takes a simpler form in the affine setting, as
illustrated in the following two examples.

\begin{exam}
\label{exam2} Consider the function
\[
f:=\sup_{t\in T}f_{t}:=\left\langle a_{t},\cdot\right\rangle -b_{t}%
,\emph{\ }a_{t}\in X^{\ast},\emph{\ }b_{t}\in\mathbb{R},
\]
where $T$ is an arbitrary set. Fix\ $x\in\operatorname*{dom}f$. By applying
(\ref{mainfb}) with $\alpha_{t}=\rho_{t,\varepsilon}=\frac{\varepsilon}%
{\max\{2f(x)-2\left\langle a_{t},x\right\rangle +2b_{t}-\varepsilon
,\varepsilon\}}$ for all $t\in T$ (see (\ref{vetb})), we obtain
\[
\mathrm{N}_{\operatorname*{dom}f}(x)=\left[  \overline{\operatorname*{co}%
}\left\{  \rho_{\varepsilon,t}a_{t}:t\in T\right\}  \right]  _{\infty},\text{
for every }\varepsilon>0,
\]
while (\ref{mainfc}) implies
\[
\mathrm{N}_{\operatorname*{dom}f}(x)=\{\theta\}\cup(\cap_{\varepsilon
>0}\overline{\operatorname*{co}}\left\{  \varepsilon\rho_{\varepsilon,t}%
a_{t}:t\in T\right\}  ).
\]
In particular, whenever $\inf_{t\in T}(\left\langle a_{t},x\right\rangle
-b_{t})>-\infty,$ Corollary \ref{corf} entails
\[
\mathrm{N}_{\operatorname*{dom}f}(x)=\left[  \overline{\operatorname*{co}%
}\left(  \cup_{t\in T}\{a_{t}\}\right)  \right]  _{\infty}=\{\theta\}\cup
(\cap_{\varepsilon>0}\overline{\operatorname*{co}}\left(  \cup_{t\in
T}\{\varepsilon a_{t}\}\right)  ).
\]

\end{exam}

\begin{exam}
Consider the support function $\sigma_{A}$ of a nonempty set $A\subset
X^{\ast}.$ Given\ $x\in\operatorname*{dom}\sigma_{A},$\ Theorem \ref{Theo13b}
yields
\begin{equation}
\mathrm{N}_{\operatorname*{dom}\sigma_{A}}(x)=\left[  \overline
{\operatorname*{co}}\left\{  \rho_{\varepsilon,t}a:a\in A\right\}  \right]
_{\infty},\text{ for every }\varepsilon>0, \label{dd}%
\end{equation}
where $\rho_{\varepsilon,t}:=\frac{\varepsilon}{\max\{2\sigma_{A}%
(x)-\left\langle 2a,x\right\rangle -\varepsilon,\varepsilon\}}.$ Moreover, if
$\sigma_{A}(-x)<+\infty,$ then
\[
\mathrm{N}_{\operatorname*{dom}\sigma_{A}}(x)=\left[  \overline
{\operatorname*{co}}(A)\right]  _{\infty}.
\]
Observe that,\ by the bipolar theorem, (\ref{dd}) is equivalent to
\[
\operatorname*{cl}\left(  \mathbb{R}_{+}(\operatorname*{dom}\sigma
_{A}-x)\right)  =\left(  \left[  \overline{\operatorname*{co}}\left\{
\rho_{\varepsilon,t}a:a\in A\right\}  \right]  _{\infty}\right)  ^{-},
\]
which, for $x=\theta,$\ simplifies to the classical
identity\ $\operatorname*{cl}\left(  \operatorname*{dom}\sigma_{A}\right)
=(\left[  \overline{\operatorname*{co}}\left(  A\right)  \right]  _{\infty
})^{-}$ (see, i.e., \cite[(3.52), p. 80]{CHLBook}).
\end{exam}

\medskip Next, we extend Theorem \ref{Theo13b} to a more general setting in
which
\[
\mathcal{P}:=\{t\in T:f_{t}\text{ is proper}\}\varsubsetneqq T,
\]
thereby\ allowing some of the (lsc convex) functions $f_{t}$ to be improper.
This extension is motivated by the fact that, in practical problems, the lower
semicontinuity assumption in Theorem \ref{Theo13b} can be replaced with the
closure condition $\bar{f}=\sup_{t\in T}\bar{f}_{t}$ (see \cite{CHLBook}).
However, in\ infinite-dimensional settings, this reformulation may produce
convex lsc functions $\bar{f}_{t}$ that\ are no longer proper, even when the
original $f_{t}$'s were proper.

The following result, whose proof relies mainly on Theorem \ref{Theo13b}%
,\ shows that improper functions remain relevant for the characterization of
$\mathrm{N}_{\operatorname*{dom}f}(x),$ contributing through their (effective) domains.

\begin{theo}
\label{Theo13} Given lsc convex\ functions $f_{t}:X\rightarrow\overline
{\mathbb{R}},\ t\in T,$ and\ $f:=\sup_{t\in T}f_{t},$ for any $x\in
\operatorname*{dom}f$ we have
\[
\mathrm{N}_{\operatorname*{dom}f}(x)=\left[  \overline{\operatorname*{co}%
}\left(  \left(  {%
{\textstyle\bigcup\limits_{t\in\mathcal{P}}}
}\partial_{\varepsilon}(\alpha_{t}f_{t})(x)\right)  {\cup}\left(  {%
{\textstyle\bigcup\limits_{t\in T\backslash\mathcal{P}}}
}\mathrm{N}_{\operatorname*{dom}f_{t}}^{\varepsilon}(x)\right)  \right)
\right]  _{\infty},
\]
for any\ $\varepsilon>0$ and $\alpha_{t}\geq\rho_{t,\varepsilon}$
$(t\in\mathcal{P})$ such that $\sup_{t\in\mathcal{P}}\alpha_{t}<+\infty$ and
$\inf_{t\in\mathcal{P}}\alpha_{t}f_{t}(x)>-\infty$ ($\rho_{t,\varepsilon,x}$
are\ defined in (\ref{vetb})).\ Consequently,
\[
\mathrm{N}_{\operatorname*{dom}f}(x)=%
{\textstyle\bigcap\limits_{\varepsilon>0}}
\overline{\operatorname*{co}}\left(  \left(  {%
{\textstyle\bigcup\limits_{t\in\mathcal{P}}}
}\partial_{\varepsilon}(\varepsilon\alpha_{t}f_{t})(x)\right)  {\cup}\left(  {%
{\textstyle\bigcup\limits_{t\in T\backslash\mathcal{P}}}
}\mathrm{N}_{\operatorname*{dom}f_{t}}^{\varepsilon}(x)\right)  \right)  ,
\]
provided that the latter set is nonempty; otherwise, $\mathrm{N}%
_{\operatorname*{dom}f}(x)=\{\theta\}.$
\end{theo}

\begin{dem}
Fix\ $x\in\operatorname*{dom}f$, $\varepsilon>0$ and let $\alpha_{t}$ be as in
the statement of the theorem. Let us first suppose that $f(x)>-\infty.$ Define
the functions $g_{t}\in\Gamma_{0}(X)$ as\ $g_{t}:=f_{t}$ if $t\in\mathcal{P},$
and $g_{t}:=\max\{f_{t},f(x)-\varepsilon\}$ if $t\in T\setminus\mathcal{P},$
and denote\ $g:=\sup_{t\in T}g_{t}.$ By the lower semi-continuity assumption,
$f$ and $g$ coincide in a neighborhood of $x$ and, so, $\mathrm{N}%
_{\operatorname*{dom}f}(x)=\mathrm{N}_{\operatorname*{dom}g}(x).$\ Moreover,
we have
\[
\frac{-\varepsilon}{2g_{t}(x)-2g(x)+\varepsilon}=\frac{-\varepsilon}%
{2f_{t}(x)-2f(x)+\varepsilon}=\rho_{t,\varepsilon},\text{ for all\ }%
t\in\mathcal{P}\setminus T_{\varepsilon}(x).
\]
In addition, since $f_{t}(x)=-\infty$ for all $t\in T\setminus\mathcal{P},$ we
have $\{t\in T:g(x)\geq g_{t}(x)-\varepsilon\}=T_{\varepsilon}(x)\cup
(T\setminus\mathcal{P}).$ Moreover, (\ref{new}) yields $\partial_{\varepsilon
}g_{t}(x)=\mathrm{N}_{\operatorname*{dom}f_{t}}^{\varepsilon}(x)$ for all
$t\in T\backslash\mathcal{P}.$ Set\ $\alpha_{t}:=1$ for $t\in T\backslash
\mathcal{P},$ so that $\inf_{t\in T}(\alpha_{t}g_{t})(x)\in\mathbb{R}$%
.$\ $Thus, applying Theorem\ \ref{Theo13b} to the $g_{t}$'s gives
\begin{align*}
\mathrm{N}_{\operatorname*{dom}f}(x) &  =\mathrm{N}_{\operatorname*{dom}%
g}(x)=\left[  \overline{\operatorname*{co}}\left(  \cup_{t\in T}%
\partial_{\varepsilon}(\alpha_{t}g_{t})(x)\right)  \right]  _{\infty}\\
&  =\left[  \overline{\operatorname*{co}}\left(  \cup_{t\in\mathcal{P}%
}\partial_{\varepsilon}(\alpha_{t}f_{t})(x){\cup}\cup_{t\in T\backslash
\mathcal{P}}\partial_{\varepsilon}g_{t}(x)\right)  \right]  _{\infty}\\
&  =\left[  \overline{\operatorname*{co}}\left(  \cup_{t\in\mathcal{P}%
}\partial_{\varepsilon}(\alpha_{t}f_{t})(x){\cup}\cup_{t\in T\backslash
\mathcal{P}}\mathrm{N}_{\operatorname*{dom}f_{t}}^{\varepsilon}(x)\right)
\right]  _{\infty},
\end{align*}
which yields the first assertion. The second formula\ then follows from the
first one by arguing as in Step 5 of the proof of Theorem \ref{Theo13b}.

Now, suppose\ that\ $f(x)=-\infty.$ Then\ $\mathcal{P}=\emptyset$ and
$f_{t}(x)=-\infty,$ for all $t\in T.$ Consider the family of lsc convex
functions $\{f_{t},$ $t\in T,$ $f_{0}\},$ where\ $f_{0}\equiv1$ (supposing
that $0\notin T$). Then $\sup_{t\in T\cup\{0\}}f_{t}(x)=1\in{\mathbb{R}}$.
Applying the previous argument to the function $h:=\sup_{t\in T\cup\{0\}}%
f_{t}$\ entails\
\begin{align*}
\mathrm{N}_{\operatorname*{dom}f}(x)  &  =\mathrm{N}_{\operatorname*{dom}%
h}(x)=\left[  \overline{\operatorname*{co}}\left(  \partial_{\varepsilon}%
f_{0}(x)\cup\left(  \cup_{t\in T}\mathrm{N}_{\operatorname*{dom}f_{t}%
}^{\varepsilon}(x)\right)  \right)  \right]  _{\infty}\\
&  =\left[  \overline{\operatorname*{co}}\left(  \cup_{t\in T}\mathrm{N}%
_{\operatorname*{dom}f_{t}}^{\varepsilon}(x)\right)  \right]  _{\infty
}=\left[  \overline{\operatorname*{co}}\left(  \cup_{t\in T\setminus
\mathcal{P}}\mathrm{N}_{\operatorname*{dom}f_{t}}^{\varepsilon}(x)\right)
\right]  _{\infty},
\end{align*}
thereby yielding the first assertion. The second formula follows analogously
to the first part of the proof.
\end{dem}

\medskip The characterization of the normal cone in Theorem \ref{Theo13b}
reveals the asymmetric role between the $\varepsilon$-active functions at $x$
and the non-$\varepsilon$-active ones; in fact, only the latter are penalized
by the weights $\rho_{t,\varepsilon}$ introduced\ in (\ref{vetb}).\ 

\begin{cor}
\label{corlaylaha} Given lsc convex functions $f_{t}:X\rightarrow
\overline{\mathbb{R}},\ t\in T,$ we define\ $f:=\sup_{t\in T}f_{t}$. Then, for
every $x\in\operatorname*{dom}f,$
\begin{align*}
\mathrm{N}_{\operatorname*{dom}f}(x)=\left[  \overline{\operatorname*{co}%
}\left( \left( \underset{t\in T_{\varepsilon}(x)}{{\cup}}\partial_{\varepsilon
}(\beta_{t}f_{t})(x)\right)\cup \left( \underset{t\in\mathcal{P}\setminus T_{\varepsilon}%
(x)}{{\cup}}\eta_{\varepsilon}\partial_{\varepsilon}(\beta_{t}\rho
_{t,\varepsilon}f_{t})(x)\right)\right.\right.
\\ 
\left.\left. \ \qquad \qquad 
\cup\left(\underset{t\in T\backslash\mathcal{P}}{{\cup}%
}\mathrm{N}_{\operatorname*{dom}f_{t}}^{\varepsilon}(x)\right)\right)  \right]
_{\infty},\text{ }%
\end{align*}
for any\ $\varepsilon>0$, any positive number $\eta_{\varepsilon}$ with
$\eta_{\varepsilon}\downarrow0,$ and any family of positive scalars
$(\beta_{t})_{t\in T}$ with $0<\inf_{t\in T}\beta_{t}\leq\sup_{t\in T}%
\beta_{t}<+\infty.$
\end{cor}

\begin{dem}
By Lemma \ref{lema1}, the weights $\rho_{t,\varepsilon}$ $(\in$ $]0,1])$
satisfy $\inf_{t\in T}\rho_{t,\varepsilon}f_{t}(x)>-\infty.$ It then\ suffices
to apply Theorem \ref{Theo13b} with $\alpha_{t}=\rho_{t,\varepsilon}$ for all
$t\in T.$ The use of the coefficients\ $\eta_{\varepsilon}$ and the $\beta
_{t}$'s is justified by (\ref{mostakim}) and Lemma \ref{ChangeParameter}, respectively.
\end{dem}

\medskip The parameters $\alpha_{t}$ appearing in Theorem\ \ref{Theo13} can be
removed by suitably modifying certain functions $f_{t}$\ without altering the
normal cone $\mathrm{N}_{\operatorname*{dom}f}(x)$. In the following
corollary, we replace\ the functions $f_{t}$ for $t\in\mathcal{P}\backslash
T_{\varepsilon}(x)$ with\
\begin{equation}
f_{t,\varepsilon}:=\max\{f_{t},f(x)-\varepsilon\}. \label{fca}%
\end{equation}
These new functions ensure that $\inf_{t\in\mathcal{P}\backslash
T_{\varepsilon}(x)}f_{t,\varepsilon}(x)>-\infty,$ while preserving the domain
of $f.$ Moreover, by \cite[Corollary 5.1.9]{CHLBook}, the set\ $\partial
_{\varepsilon}f_{t,\varepsilon}(x)\ $can be expressed by means of\ $f_{t}$ as%
\[
\partial_{\varepsilon}f_{t,\varepsilon}(x)={%
{\textstyle\bigcup\limits_{\lambda\in]0,1]}}
}\lambda\partial_{\varepsilon_{\lambda}}f_{t}(x)\cup\mathrm{N}%
_{\operatorname*{dom}f_{t}}^{\varepsilon}(x),
\]
where $\varepsilon_{\lambda}:=\lambda^{-1}\varepsilon+f_{t}%
(x)-f(x)+\varepsilon.$

\begin{cor}
\label{corThm13a} Let\ functions $f_{t}:X\rightarrow\overline{\mathbb{R}%
},\ t\in T,$ be convex and lsc, and denote\ $f:=\sup_{t\in T}f_{t}.$ Let
$x\in\operatorname*{dom}f.$ Then, for every\ $\varepsilon>0,$ we have
\begin{align*}
\mathrm{N}_{\operatorname*{dom}f}(x)=\left[  \overline{\operatorname*{co}%
}\left(  \left(  \underset{t\in T_{\varepsilon}(x)}{{\cup}}\partial
_{\varepsilon}f_{t}(x)\right)  {\cup}\left(  \underset{t\in\mathcal{P}%
\backslash T_{\varepsilon}(x)}{{\cup}}\partial_{\varepsilon}f_{t,\varepsilon
}(x)\right) 
\right.\right.
\\ 
\left.\left.
{\cup}\left(  \underset{t\in T\backslash\mathcal{P}}{{\cup}%
}\mathrm{N}_{\operatorname*{dom}f_{t}}^{\varepsilon}(x)\right)  \right)
\right]  _{\infty}. \label{bb}%
\end{align*}

\end{cor}

\begin{dem}
Fix $x\in\operatorname*{dom}f$ and $\varepsilon>0.$ First, we assume that
$f(x)\in\mathbb{R}$. We consider\ the convex functions $g,$ $g_{t}%
:X\rightarrow\overline{\mathbb{R}}$, $t\in T,$ defined as\
\[
g_{t}:=f_{t}\text{ if }t\in T_{\varepsilon}(x)\cup(T\setminus\mathcal{P}%
),\text{ }g_{t}:=f_{t,\varepsilon}\text{ if }t\in\mathcal{P}\backslash
T_{\varepsilon}(x),\text{ }g:=\sup\nolimits_{t\in T}g_{t}.
\]
Observe that $g_{t}\in\Gamma_{0}(X)$ for all $t\in\mathcal{P},$ and
$\inf_{t\in\mathcal{P}}g_{t}(x)\geq f(x)-\varepsilon>-\infty.$ Moreover,
because $f$ is lsc and $f(x)\in\mathbb{R},$ the functions $f$ and $g$ coincide
in some neighborhood of $x,$ implying that\ $\mathrm{N}_{\operatorname*{dom}%
f}(x)=\mathrm{N}_{\operatorname*{dom}g}(x).$ Therefore the conclusion follows by
applying Theorem \ref{Theo13} to the family $\{g_{t},$ $t\in T\}$ with
$\alpha_{t}=1,$ $t\in\mathcal{P}.$

If $f(x)=-\infty,$ then the conclusion follows by proceeding as in the proof
of Theorem \ref{Theo13}, taking the family $\{f_{t},$ $t\in T,$ $f_{0}\}$
with\ $f_{0}\equiv1.$
\end{dem}

\medskip

Further characterizations of the normal cone are presented in Section
\ref{secultim2}, offering refined formulas for the normal cone when\ the
supremum function is continuous at some point.

\section{Subdifferential of the supremum
function\label{SectionSubdifferential}}

The results of the previous section on normal cones are applied here to
characterize the subdifferential set $\partial f(x)$ of the supremum function
\[
f:=\sup_{t\in T}f_{t},
\]
where\ $\{f_{t}:X\rightarrow\overline{\mathbb{R}},$ $t\in T\}$ is an arbitrary
family of convex functions defined on a lcs $X$. As is evident from the known
characterizations of $\partial f(x)$ (see, e.g., \cite{Br72, HLZ08, Ioffe12,
LoVo10} and the recent book \cite{CHLBook}), a dominant role in the
construction of $\partial f(x)$ is played by the functions $f_{t}$ indexed by
the set of $\varepsilon$-active indices at $x\in\operatorname*{dom}f,$
\[
T_{\varepsilon}(x):=\{t\in\mathcal{P}:f_{t}(x)\geq f(x)-\varepsilon\},\text{
}\varepsilon>0,
\]
where $\mathcal{P}$ denotes the collection of proper functions.$\ $Our
analysis further shows that the remaining proper functions also contribute
through a suitable penalization given by the weights $\rho_{t,\varepsilon}$
defined in (\ref{vetb}),
\[
\rho_{t,\varepsilon}:=\frac{-\varepsilon}{2f_{t}(x)-2f(x)+\varepsilon},\text{
}t\in\mathcal{P}\setminus T_{\varepsilon}(x)\text{ and }\rho_{t,\varepsilon
}:=1\text{ for }t\in T_{\varepsilon}(x).
\]
Improper functions, on the other hand, enter the characterization of $\partial
f(x)$ through their corresponding $\varepsilon$-normal sets. The following
theorem provides a full characterization of the subdifferential of the
supremum function $f,$ highlighting the following components:%
\begin{equation}
A_{\varepsilon}(x):=%
{\textstyle\bigcup\limits_{t\in T_{\varepsilon}(x)}}
\partial_{\varepsilon}f_{t}(x),\text{ }B_{\varepsilon,\alpha}(x):=%
{\textstyle\bigcup\limits_{t\in\mathcal{P}\backslash T_{\varepsilon}(x)}}
\partial_{\varepsilon}(\alpha_{t}f_{t})(x),\text{ }C_{\varepsilon}(x):=%
{\textstyle\bigcup\limits_{t\in T\backslash\mathcal{P}}}
\mathrm{N}_{\operatorname{dom}f_{t}}^{\varepsilon}(x).\label{akbal}%
\end{equation}
Notice that the weights $\rho_{t,\varepsilon}$ provide an example of the
parameters $\alpha_{t}$ used below; see Lemma \ref{lema1}.

\begin{theo}
\label{thmsub}\ Let $f_{t}:X\rightarrow\overline{\mathbb{R}},$ $t\in T,$ be
convex and lsc, and denote $f:=\sup_{t\in T}f_{t}.$ Then, for every
$x\in\operatorname{dom}f,$\ we have
\begin{equation}
\partial f(x)=\bigcap_{\varepsilon>0}\overline{\operatorname*{co}}\left(
A_{\varepsilon}(x)+\varepsilon(B_{\varepsilon,\alpha}(x)\cup C_{\varepsilon
}(x)\cup\{\theta\})\right)  ,\label{FGSubSup}%
\end{equation}
for any $\alpha_{t}\geq\rho_{t,\varepsilon}$ ($t\in\mathcal{P}\backslash
T_{\varepsilon}(x)$) such that $\sup\limits_{t\in\mathcal{P}\backslash
T_{\varepsilon}(x)}\alpha_{t}<+\infty$ and $\inf\limits_{t\in\mathcal{P}%
\backslash T_{\varepsilon}(x)}\alpha_{t}f_{t}(x)>-\infty$.
\end{theo}

\begin{rem}
[before the proof]As shown in the proof, the inclusion \textquotedblleft%
$\subset$\textquotedblright\ in (\ref{FGSubSup}) holds without requiring that
$\inf_{t\in T}\alpha_{t}f_{t}(x)>-\infty;$ the latter being necessary\ only
for the opposite inclusion. Furthermore, by Lemma \ref{ChangeParameter},
(\ref{FGSubSup}) continues to hold if, instead of $\alpha_{t}\ $%
($t\in\mathcal{P}\backslash T_{\varepsilon}(x)$), we use any scalar multiple
$\beta_{t}\alpha_{t}\ $with $0<\inf_{t\in\mathcal{P}\backslash T_{\varepsilon
}(x)}\beta_{t}$ and $\sup_{t\in\mathcal{P}\backslash T_{\varepsilon}(x)}%
\beta_{t}<+\infty.$
\end{rem}

\begin{dem}
Fix $x\in\operatorname{dom}f.$ We may assume throughout the proof that
$f(x)\in\mathbb{R};$ otherwise, we would have $f_{t}(x)=f(x)=-\infty$ for all
$t\in T,$ and all the sets $\partial f(x)$ and $A_{\varepsilon}(x),$
$\varepsilon>0,$ would be\ empty. In that case, (\ref{FGSubSup}) holds
trivially. Moreover, up to replacing each $f_{t}$ with $f_{t}-f(x),$ we may
also suppose, without loss of generality, that $f(x)=0.$

Let us first prove the inclusion \textquotedblleft$\subset$\textquotedblright%
\ of (\ref{FGSubSup}) in the\ nontrivial case where $\partial f(x)\not =%
\emptyset.$ Fix $\varepsilon>0$ and let $\alpha$ be as in the statement of the
theorem. Pick a $\theta$-neighborhood\ $U\in\mathcal{N}$ and\ choose
$L\in\mathcal{F}(x)$ such that $L^{\perp}\subset U.$ Then, by applying Theorem
\ref{Theo13} to the family $\{f_{t},$ $t\in T;$ $\mathrm{I}_{L}\},$ we obtain
\[
\mathrm{N}_{L\cap\operatorname*{dom}f}(x)\subset\left[  \overline
{\operatorname*{co}}\left(  A_{\varepsilon}(x)\cup L^{\perp}\cup\left(
B_{\varepsilon,\alpha}(x)\cup C_{\varepsilon}(x)\cup\{\theta\}\right)
\right)  \right]  _{\infty},
\]
or equivalently, thanks to\ (\ref{mostakim}),
\begin{equation}
\mathrm{N}_{L\cap\operatorname*{dom}f}(x)\subset\left[  \overline
{\operatorname*{co}}\left(  A_{\varepsilon}(x)+L^{\perp}+\varepsilon\left(
B_{\varepsilon,\alpha}(x)\cup C_{\varepsilon}(x)\cup\{\theta\}\right)
\right)  \right]  _{\infty}.\label{est}%
\end{equation}
Observe that $T_{\varepsilon}(x)\neq\emptyset$ and $A_{\varepsilon}%
(x)\neq\emptyset,$ since $f_{t}\in\Gamma_{0}(X)\ $for all $t\in T_{\varepsilon
}(x).$ Therefore, combining (\ref{est}) and (\ref{refdem16}),
\begin{align*}
\partial f(x) &  \subset\overline{\operatorname*{co}}\left(  A_{\varepsilon
}(x)+\mathrm{N}_{L\cap\operatorname*{dom}f}(x)\right)  \\
&  \subset\overline{\operatorname*{co}}\left(  A_{\varepsilon}(x)+\left[
\overline{\operatorname*{co}}\left(  A_{\varepsilon}(x)+\varepsilon
(B_{\varepsilon,\alpha}(x)\cup C_{\varepsilon}(x)\cup\{\theta\})+L^{\perp
}\right)  \right]  _{\infty}\right)  \\
&  \subset\overline{\operatorname*{co}}\left(  A_{\varepsilon}(x)+\varepsilon
(B_{\varepsilon,\alpha}(x)\cup C_{\varepsilon}(x)\cup\{\theta\})+L^{\perp
}\ \right)  \\
&  \subset\operatorname*{co}\left(  A_{\varepsilon}(x)+\varepsilon
(B_{\varepsilon,\alpha}(x)\cup C_{\varepsilon}(x)\cup\{\theta\})\right)  +2U,
\end{align*}
and \textquotedblleft$\subset$\textquotedblright\ in (\ref{FGSubSup}) follows
by intersecting\ first over $U\in\mathcal{N}$ and next over $\varepsilon>0.$

To establish the opposite inclusion, we fix\ $y\in\operatorname*{dom}f$
$(\subset\cap_{t\in T}\operatorname*{dom}f_{t})$ and $\varepsilon>0.$ Then
$\mathrm{\sigma}_{C_{\varepsilon}(x)}(y-x)\leq\varepsilon,$ $\mathrm{\sigma
}_{A_{\varepsilon}(x)}(y-x)\leq\sup\nolimits_{t\in T_{\varepsilon}(x)}%
(f_{t}(y)-f_{t}(x)+\varepsilon)\leq f(y)-f(x)+2\varepsilon,$ and
\[
\mathrm{\sigma}_{B_{\varepsilon,\alpha}(x)}(y-x)\leq\sup\nolimits_{t\in
\mathcal{P}\setminus T_{\varepsilon}(x)}(\alpha_{t}f_{t}(y)-\alpha_{t}%
f_{t}(x)+\varepsilon)\leq(M_{1}f^{+}(y)-M_{2}\alpha_{t}f_{t}(x)+\varepsilon
)^{+},
\]
where $M_{1}:=\sup_{t\in\mathcal{P}\setminus T_{\varepsilon}(x)}\alpha_{t}$
and $M_{2}:=\inf_{t\in\mathcal{P}\setminus T_{\varepsilon}(x)}\alpha_{t}%
f_{t}(x).$ Finally,\ if\ $E$ is the right hand side of (\ref{FGSubSup}), then
we obtain
\[
\mathrm{\sigma}_{E}(y-x)\leq f(y)-f(x)+\max\{\varepsilon(M_{1}f^{+}%
(y)-M_{2}f_{t}(x)+\varepsilon)^{+},\varepsilon^{2}\}+2\varepsilon,
\]
and the inclusion \textquotedblleft$\supset$\textquotedblright\ in
(\ref{FGSubSup}) follows when $\varepsilon\downarrow0.$
\end{dem}

\medskip We give an extension of Brøndsted's formula, originally established
for finite collections in \cite{Br72}. A finite-dimensional variant of
(\ref{t2}) was also proven in\ \cite[Proposition 6.3]{HL08} under the
continuity of the $f_{t}$'s. Instead of the lower semicontinuity requirement
in Theorem \ref{Theo13}, we use here the following closure condition $\bar
{f}=\sup_{t\in T}\bar{f}_{t}.$ This\ naturally holds when the $f_{t}$'s are
lsc, but it is also satisfied in other notable cases, such as when the
supremum function $f$ is continuous at some point\ (see \cite[Corollary
9]{HLZ08}, \cite{Ro70}).

\begin{cor}
\label{bronds}Given convex\ functions $f_{t},$ $t\in T,$ and$\ f:=\sup_{t\in
T}f_{t},$ we assume that $\bar{f}=\sup_{t\in T}\bar{f}_{t}$.\ If $x\in X$ is
such that $\bar{f}_{t}(x)=f(x)\in\mathbb{R}$, for all $t\in T,$ then
\begin{equation}
\partial f(x)=\bigcap_{\varepsilon>0}\overline{\operatorname*{co}}\left(  {%
{\textstyle\bigcup\limits_{t\in T}}
}\partial_{\varepsilon}f_{t}(x)\right)  . \label{t2}%
\end{equation}

\end{cor}

\begin{dem}
Firstly, we have $\partial_{\varepsilon}f_{t}(x)\subset\partial_{\varepsilon
}f(x)$, for all $t\in T$ and $\varepsilon>0,$ and so
\[
{\cap}_{\varepsilon>0}\overline{\operatorname*{co}}\left(  {\cup}_{t\in
T}\partial_{\varepsilon}f_{t}(x)\right)  \subset{\cap}_{\varepsilon>0}%
\partial_{\varepsilon}f(x)=\partial f(x).
\]
Thus, the inclusion \textquotedblleft$\supset$\textquotedblright\ in
(\ref{t2})\ holds. To prove the opposite inclusion we may assume that
$\partial f(x)\neq\emptyset;$ otherwise, (\ref{t2}) trivially holds. This
implies that $\partial f(x)=\partial\bar{f}(x).$ Since $\bar{f}_{t}%
(x)=f(x)\geq f_{t}(x),$ it follows that $\partial_{\varepsilon}\bar{f}%
_{t}(x)=\partial_{\varepsilon}f_{t}(x)$, for all $t\in T$ and $\varepsilon>0.$
Therefore the conclusion follows by applying Theorem \ref{thmsub} to the
family of proper, lsc, and convex functions $\{\bar{f}_{t},$ $t\in T\},$%
\[
\partial f(x)=\partial\bar{f}(x)={\cap}_{\varepsilon>0}\overline
{\operatorname*{co}}\left(  {\cup}_{t\in T}\partial_{\varepsilon}\bar{f}%
_{t}(x)+\{\theta\}\right)  ={\cap}_{\varepsilon>0}\overline{\operatorname*{co}%
}\left(  {\cup}_{t\in T}\partial_{\varepsilon}f_{t}(x)\right)  .
\]

\end{dem}

\medskip We give an illustration\ of Theorem \ref{thmsub} in\ the affine case.
The case where the supremum function is continuous at some point will be
examined in Example \ref{exammcont}.

\begin{exam}
\label{examm}Given a nonempty set $T,$ we consider the function
\[
f:=\sup_{t\in T}\left\{  \left\langle a_{t},\cdot\right\rangle -b_{t}\right\}
,\text{ }a_{t}\in X^{\ast},\text{ }b_{t}\in\mathbb{R}.
\]
Then, combining Theorem \ref{thmsub} and Lemma \ref{lema1},\ for all $x\in
f^{-1}(0)$ we obtain
\[
\partial f(x)=\cap_{\varepsilon>0}\overline{\operatorname*{co}}(\left\{
a_{t}:\left\langle a_{t},x\right\rangle -b_{t}\geq-\varepsilon\right\}
+\varepsilon(B_{\varepsilon}(x)\cup\{0\})),\text{ }%
\]
where $B_{\varepsilon}(x):=\left\{  \frac{\varepsilon a_{t}}{2b_{t}%
-2\left\langle a_{t},x\right\rangle -\varepsilon}:\left\langle a_{t}%
,x\right\rangle -b_{t}<-\varepsilon\right\}  .$
\end{exam}

Further refinements are presented in Section \ref{secultim2}:
Theorem\ \ref{conecont} addresses the case where the supremum function $f$ is
continuous at some point, while Corollary\ \ref{olab} concerns the so-called
compact-continuous setting.

\section{Optimality conditions\label{secopt}}

In this section, we derive general\ optimality conditions for the following
convex infinite convex program, posed in a lcs $X$,\ %

\[%
\begin{array}
[c]{llll}%
(\mathtt{P}) & \inf & f_{0}(x) & \text{s. t. }f_{t}(x)\leq0,\text{ }t\in T,
\end{array}
\]
where each $f_{t}:X\rightarrow\overline{\mathbb{R}}$, $t\in T\cup\{0\}$
($T\neq\emptyset,$ $0\notin T$), is supposed to be convex and lsc. We suppose
that $(\mathtt{P})$ has a finite optimal value $v(\mathtt{P})$, and we
denote\
\[
f:=\sup_{t\in T}f_{t}.
\]
The analysis of optimality conditions for\ $(\mathtt{P})$ relies on how the
constraints are penalized. In this paper, we base on the so-called performance
function $\max\{f_{0}-v(\mathtt{P}),$ $f_{t},$ $t\in T\}.$ Next, we present
Karush-Kuhn-Tucker (KKT) type optimality conditions for problem $(\mathtt{P})$
under the Slater condition\
\[
\lbrack f:=\sup_{t\in T}f_{t}<0]\cap\operatorname*{dom}f_{0}\neq\emptyset.
\]
The condition below is formulated explicitly in terms of $\varepsilon$-active,
non-$\varepsilon$-active, and improper constraints, respectively represented
for $\varepsilon>0$ by
\[
A_{\varepsilon}(x):=\underset{t\in T_{\varepsilon}(x)}{\cup}\partial
_{\varepsilon}f_{t}(x),\text{ }\mathcal{B}_{\varepsilon}(x):=\underset{t\in
T\backslash T_{\varepsilon}(x)}{\cup}\varepsilon\partial_{\varepsilon}%
(\rho_{t,\varepsilon}f_{t})(x),
\]
and $$C_{\varepsilon}(x):=\underset{t\in T\backslash\mathcal{P}}{\cup
}\mathrm{N}_{\operatorname*{dom}f_{t}}^{\varepsilon}(x),$$ where $\mathcal{P}%
=\{t\in T:f_{t}$ is proper$\},$ $T_{\varepsilon}(x)=\{t\in \mathcal{P}:f_{t}(x)\geq
f(x)-\varepsilon\}$ and $\rho_{t,\varepsilon}$ is defined as in (\ref{vetb}). The theorem
then extends the corresponding result in \cite[Theorem 8.2.6 ]{CHLBook},
originally proved in the compact-continuous setting.

\begin{theo}
\label{thmoptb}Let\ $f_{0}$ and $f_{t},\ t\in T,$ be convex and lsc functions
such that $v(\mathtt{P})\in\mathbb{R}$, and denote $f:=\sup_{t\in T}f_{t}.$
Assume\ that\emph{:}

$\emph{(i)}$ either $f_{0}$ or $f$ is continuous at some point in the domain
of the other, and

$\emph{(ii)}$ the Slater condition holds.\newline Then a feasible point $x$ is
optimal for\ $(\mathtt{P})$ if and only if:\newline$(a)$ when $f(x)<0,$
\begin{equation}
\label{th1}0\in\partial f_{0}(x)+\left[  \overline{\operatorname*{co}}\left(
A_{\varepsilon}(x)\cup\mathcal{B}_{\varepsilon}(x)\cup C_{\varepsilon
}(x)\right)  \right]  _{\infty},\text{ for all }\varepsilon>0,
\end{equation}
$(b)$ when $f(x)=0,$ either (\ref{th1}) holds, or the exists some $\lambda>0$
such that
\begin{equation}
\label{th2}0\in\lambda\partial f_{0}(x)+%
{\textstyle\bigcap\limits_{\varepsilon>0}}
\overline{\operatorname*{co}}\left(  A_{\varepsilon}(x)+\varepsilon
(\mathcal{B}_{\varepsilon}(x)\cup C_{\varepsilon}(x)\cup\{\theta\})\right)  .
\end{equation}

\end{theo}

\begin{dem}
Assume that $x$ is optimal for\ $(\mathtt{P}),$ and set $g:=\max
\{f_{0}-v(\mathtt{P}),f\}.$ Then $g(x)=0$\ and $\theta\in\partial g(x).$
When\ $f(x)<0,$ \cite[Corollary 6.5.3]{CHLBook} implies $\theta\in\partial
g(x)=\partial f_{0}(x)+\mathrm{N}_{\operatorname{dom}f}(x)$ and (\ref{th1})
follows by Theorem \ref{Theo13b} with $\alpha_{t}=\rho_{t,\varepsilon}$.
When\ $f(x)=0,$ \cite[Corollary 6.5.3]{CHLBook} yields%
\[
\theta\in\partial g(x)=\operatorname*{co}\{\partial f_{0}(x),\partial
f(x)\}+\mathrm{N}_{\operatorname{dom}f_{0}}(x)+\mathrm{N}_{\operatorname{dom}%
f}(x),
\]
and three cases my occur: If\ $\theta\in\partial f_{0}(x)+\mathrm{N}%
_{\operatorname{dom}f}(x),$ then we are in the situation above where
(\ref{th1}) holds. If $\theta\in\partial f(x)+\mathrm{N}_{\operatorname{dom}%
f_{0}}(x),$ then we get the contradiction $0=f(x)\leq f(x_{0})<0$ for every
Slater point $x_{0}\in\operatorname*{dom}f_{0}.$ If there exists
some\ $\lambda>0$ such that $\theta\in\lambda\partial f_{0}(x)+\partial f(x),$
then (\ref{th2}) follows by Theorem \ref{thmsub}.

The proof of the sufficiency of the given condition follows by the same
arguments as in Step 4 of\ the proof of Theorem \ref{Theo13b}, while also
making use of (\ref{4.5}) from\ Lemma \ref{lema1}.
\end{dem}

\medskip We apply Theorem \ref{thmoptb} to the particular case of an infinite
linear optimization problem:%
\[%
\begin{array}
[c]{llll}%
(\mathtt{PL}) & \inf & \left\langle c,x\right\rangle  & \text{s. t.
}\left\langle a_{t},x\right\rangle \leq b_{t},\text{ }t\in T,
\end{array}
\]
where $c,a_{t}\in X^{\ast},$ $t\in T.$ Let $x\in X$ be a feasible point $x$
for\ $(\mathtt{PL}).$ Under the Slater condition, Theorem \ref{thmoptb} and
Example \ref{examm} entail the following result, where we denote
$f:=\sup_{t\in T}(\left\langle a_{t},\cdot\right\rangle -b_{t})$ and
$T_{\varepsilon}(x):=\{t\in T:\left\langle a_{t},x\right\rangle -b_{t}\geq
f(x)-\varepsilon\}.$

\begin{cor}
The point\ $x$ is optimal for\ $(\mathtt{PL})$ if and only if:\newline$(a)$
when $f(x)<0,$
\begin{equation}
-c\in\left[  \overline{\operatorname*{co}}\left(  \frac{\varepsilon a_{t}%
}{\max\{2f(x)-2(\left\langle a_{t},x\right\rangle -b_{t})-\varepsilon
,\varepsilon\}}:t\in T\right)  \right]  _{\infty},\text{ for all }%
\varepsilon>0, \label{th1a}%
\end{equation}
$(b)$ when $f(x)=0,$ either (\ref{th1a}) holds, or the exists some $\lambda>0$
such that
\[
-\lambda c\in%
{\textstyle\bigcap\limits_{\varepsilon>0}}
\overline{\operatorname*{co}}\left(  \left\{  a_{t}:t\in T_{\varepsilon
}(x)\right\}  +\left\{  \frac{-\varepsilon^{2}a_{t}}{2(\left\langle
a_{t},x\right\rangle -b_{t})-2f(x)+\varepsilon}\cup\{\theta\}:t\notin 
T_{\varepsilon}(x)\right\}  \right)  .
\]

\end{cor}

\section{Enhanced structural framework\label{secultim2}}

This section gathers several refinements of the previous characterizations,
assuming\ that the supremum function $$f:=\sup_{t\in T}f_{t}$$ is continuous at
some point, or that the underlying framework is the so-called
compact-continuous setting. We assume that the $f_{t}$'s are convex and lsc
functions defined on the lcs $X$.

We start by revisiting the characterization of the normal cone under the
continuity of the supremum function $f:=\sup_{t\in T}f_{t}$ at some point,
which does not need to be the nominal one. In this case, the normal cone
admits a simpler representation, with a clear separation into three components
that highlights the asymmetric role of $\varepsilon$-active and
non-$\varepsilon$-active functions at the nominal point $x\in
\operatorname*{dom}f,$ through the sets
\begin{equation}
A_{\varepsilon}(x):=\underset{t\in T_{\varepsilon}(x)}{\cup}\partial
_{\varepsilon}f_{t}(x),\text{ }B_{\varepsilon,\alpha}(x):=\underset
{t\in\mathcal{P}\setminus T_{\varepsilon}(x)}{\cup}\partial_{\varepsilon
}(\alpha_{t}f_{t})(x),\text{ }\varepsilon>0, \label{sets}%
\end{equation}
where $\alpha=\alpha(\varepsilon)\in\mathbb{R}_{+}^{\mathcal{P}\setminus
T_{\varepsilon}(x)},$ $\mathcal{P}=\{t\in T:f_{t}$ is proper$\}$ and
$T_{\varepsilon}(x)=\{t\in \mathcal{P}:f_{t}(x)\geq f(x)-\varepsilon\}.$ As before,
improper functions are still present in our characterizations, appearing in
the sets
\begin{equation}
C_{\varepsilon}(x):=\underset{t\in T\setminus\mathcal{P}}{\cup}\mathrm{N}%
_{\operatorname*{dom}f_{t}}^{\varepsilon}(x),\text{ }\varepsilon>0.
\label{setsb}%
\end{equation}
Recall that the weights $\rho_{\varepsilon}\in\lbrack0,1]^{\mathcal{P}}$ are defined by
$\rho_{t,\varepsilon}:=\varepsilon(\max\{2f(x)-2f_{t}(x)-\varepsilon
,\varepsilon\})^{-1},$ $t\in \mathcal{P}$ (see (\ref{vetb})).

\begin{theo}
\label{conecont} Let $f_{t}:X\rightarrow\overline{\mathbb{R}},\ t\in T,$
be\ lsc convex\ functions such that $f:=\sup_{t\in T}f_{t}$ is continuous
somewhere in $\operatorname*{dom}f.$ Then, for every $x\in f^{-1}%
(\mathbb{R}),$\ we have
\begin{equation}
\mathrm{N}_{\operatorname*{dom}f}(x)=[A_{\varepsilon}(x)]_{\infty
}+[B_{\varepsilon,\alpha}(x)]_{\infty}+[C_{\varepsilon}(x)]_{\infty},
\label{f1}%
\end{equation}
for any\ $\varepsilon>0$ and any $\alpha=\alpha(\varepsilon)\in\mathbb{R}%
_{+}^{\mathcal{P}\setminus T_{\varepsilon}(x)}$ such that $\alpha_{t}\geq
\rho_{t,\varepsilon}$, $\sup_{t\in\mathcal{P}\setminus T_{\varepsilon
}(x)}\alpha_{t}<+\infty$\ and $\inf_{t\in\mathcal{P}\setminus T_{\varepsilon}%
(x)}\alpha_{t}f_{t}(x)>-\infty$. Consequently,
\begin{equation}
\mathrm{N}_{\operatorname*{dom}f}(x)=\left[  \underset{\varepsilon>0}{\cap
}\overline{\operatorname*{co}}(A_{\varepsilon}(x))\right]  _{\infty}+\left(
\underset{\delta>0}{\cap}\operatorname*{cl}\left(  \underset{0<\varepsilon
<\delta}{\cup}[\overline{\operatorname*{co}}(B_{\varepsilon,\alpha
}(x))]_{\infty}\right)  \right)  +\{\theta\}\cup\underset{\varepsilon>0}{\cap
}\overline{\operatorname*{co}}(C_{\varepsilon}(x)), \label{f2}%
\end{equation}
with $\alpha=\alpha(\varepsilon)\in\mathbb{R}_{+}^{\mathcal{P}\setminus
T_{\varepsilon}(x)}$ being as above.
\end{theo}

\begin{dem}
Fix $x\in f^{-1}(\mathbb{R})$ and $\varepsilon>0$.$\ $We define\ the convex
functions%
\[
f_{\varepsilon}:=\sup_{t\in\mathcal{P}\setminus T_{\varepsilon}(x)}%
f_{t},\text{ }g_{\varepsilon}:=\sup_{t\in T_{\varepsilon}(x)}f_{t},\text{ and
}h:=\sup_{t\in T\setminus\mathcal{P}}f_{t},
\]
and, for notational convenience, we denote
\[
A_{\varepsilon}:=\overline{\operatorname*{co}}(A_{\varepsilon}(x)),\text{
}C_{\varepsilon}:=\overline{\operatorname*{co}}(C_{\varepsilon}(x)),\text{ and
}B_{\varepsilon,\alpha}:=\overline{\operatorname*{co}}(B_{\varepsilon,\alpha
}(x))\text{ for }\alpha\in\mathbb{R}_{+}^{\mathcal{P}\setminus T_{\varepsilon
}(x)}.
\]
Since\ $\operatorname*{dom}f=\operatorname*{dom}f_{\varepsilon}\cap
\operatorname*{dom}g_{\varepsilon}\cap\operatorname*{dom}h$ and $f$ is
continuous a some point, applying the Moreau-Rockafellar Theorem (e.g.,
\cite[Proposition 4.1.16]{CHLBook}) together with Theorem \ref{Theo13b}
yields\
\begin{equation}
\mathrm{N}_{\operatorname*{dom}f}(x)=\mathrm{N}_{\operatorname*{dom}%
f_{\varepsilon}}(x)+\mathrm{N}_{\operatorname*{dom}g_{\varepsilon}%
}(x)+\mathrm{N}_{\operatorname*{dom}h}(x)=[B_{\varepsilon,\tilde{\rho
}_{\varepsilon}}]_{\infty}+[A_{\varepsilon}]_{\infty}+[C_{\varepsilon
}(x)]_{\infty}, \label{s1}%
\end{equation}
where $\tilde{\rho}_{\varepsilon}:=(\tilde{\rho}_{t,\varepsilon}%
)\in~]0,1]^{\mathcal{P}\setminus T_{\varepsilon}(x)}$ defines\ the weights
associated with the supremum function $f_{\varepsilon},$ given by
\begin{equation}
\tilde{\rho}_{t,\varepsilon}:=\varepsilon(\max\{2f_{\varepsilon}%
(x)-2f_{t}(x)-\varepsilon,\varepsilon\})^{-1},\text{ }t\in\mathcal{P}\setminus
T_{\varepsilon}(x). \label{rt}%
\end{equation}
Let $\alpha\in\mathbb{R}_{+}^{\mathcal{P}\setminus T_{\varepsilon}(x)}$ be as
in the statement of the theorem. If $\mathcal{P}\setminus T_{\varepsilon
}(x)=\emptyset,$ then $[B_{\varepsilon,\tilde{\rho}_{\varepsilon}}]_{\infty
}=[B_{\varepsilon,\alpha}]_{\infty}=\{\theta\}$ and (\ref{f1}) follows from
(\ref{s1}). Otherwise, $\mathcal{P}\setminus T_{\varepsilon}(x)\neq\emptyset$
and Lemma \ref{lema2}$(i)$ yields\ some $M\geq1$ such that $\rho
_{t,\varepsilon}\leq\tilde{\rho}_{t,\varepsilon}\leq M\rho_{t,\varepsilon}\leq
M\alpha_{t}$ for all $t\in\mathcal{P}\setminus T_{\varepsilon}(x)$. Based on
this, Lemma \ref{ChangeParameter} and Theorem \ref{Theo13b} entail\
\begin{equation}
\mathrm{N}_{\operatorname*{dom}f_{\varepsilon}}(x)=[B_{\varepsilon,\tilde
{\rho}_{\varepsilon}}]_{\infty}=[B_{\varepsilon,\rho_{\varepsilon}}]_{\infty
}=[B_{\varepsilon,\alpha}(x)]_{\infty}, \label{lc}%
\end{equation}
and (\ref{f1}) follows by (\ref{s1}).

\medskip To prove\ (\ref{f2}), we start from relation (\ref{s1}),
\begin{equation}
\mathrm{N}_{\operatorname*{dom}f}(x)=[A_{\varepsilon}]_{\infty}%
+[C_{\varepsilon}]_{\infty}+[B_{\varepsilon,\tilde{\rho}_{\varepsilon}%
}]_{\infty},\text{ for all }\varepsilon>0, \label{abo}%
\end{equation}
where $\tilde{\rho}_{\varepsilon}:=(\tilde{\rho}_{t,\varepsilon}%
)_{t\in\mathcal{P}\setminus T_{\varepsilon}(x)}$ is defined in (\ref{rt}). It
is clear that the families $([A_{\varepsilon}]_{\infty})_{\varepsilon>0}$ and
$([C_{\varepsilon}]_{\infty})_{\varepsilon>0}$ are non-increasing as
$\varepsilon\downarrow0.$ At the same time, Lemma \ref{lema2}$(ii)$ ensures
that the family $([B_{\varepsilon,\tilde{\rho}_{\varepsilon}}]_{\infty
})_{\varepsilon>0}$ is non-decreasing\ as $\varepsilon\downarrow0$. Hence, by
applying the support function in (\ref{abo}) and, next,\ taking\ the upper
limits as $\varepsilon\downarrow0$ we get
\begin{equation}
\sigma_{\mathrm{N}_{\operatorname*{dom}f}(x)}=\inf_{\varepsilon>0}%
\sigma_{\lbrack A_{\varepsilon}]_{\infty}}+\inf_{\varepsilon>0}\sigma_{\lbrack
C_{\varepsilon}]_{\infty}}+\inf_{\delta>0\text{ }}\sup_{0<\varepsilon<\delta
}\sigma_{\lbrack B_{\varepsilon,\tilde{\rho}_{\varepsilon}}]_{\infty}}.
\label{yy}%
\end{equation}
Since $\operatorname*{dom}f\subset\lbrack\sigma_{\mathrm{N}%
_{\operatorname*{dom}f}(x)}(\cdot-x)\leq0]$ and $f$ is continuous on the
interior of its domain, it follows that the (convex) function $\sigma
_{\mathrm{N}_{\operatorname*{dom}f}(x)}(\cdot-x)$ is also continuous on
$\operatorname*{int}(\operatorname*{dom}f).$ Moreover, since the three infimum
functions appearing on the right-hand side of (\ref{yy}) are convex and
nonnegative, they are dominated\ by $\sigma_{\mathrm{N}_{\operatorname*{dom}%
f}(x)}$ and, so, are also continuous on $\operatorname*{int}%
(\operatorname*{dom}f)-x$ $(\neq\emptyset).$ Therefore, applying
\cite[Proposition 2.2.11]{CHLBook}, (\ref{yy}) yields
\[
\sigma_{\mathrm{N}_{\operatorname*{dom}f}(x)}=\operatorname*{cl}%
(\inf_{\varepsilon>0}\sigma_{\lbrack A_{\varepsilon}]_{\infty}}%
)+\operatorname*{cl}(\inf_{\varepsilon>0}\sigma_{\lbrack C_{\varepsilon
}]_{\infty}})+\operatorname*{cl}(\inf_{\delta>0}\sigma_{\cup_{0<\varepsilon
<\delta}[B_{\varepsilon,\tilde{\rho}_{\varepsilon}}]_{\infty}}).
\]
Moreover, since the last three infimum functions are also proper
(nonnegative), by \cite[Proposition 3.2.8(ii)]{CHLBook} we conclude that
\begin{align}
\sigma_{\mathrm{N}_{\operatorname*{dom}f}(x)}  &  =\sigma_{\cap_{\varepsilon
>0}[A_{\varepsilon}]_{\infty}}+\sigma_{\cap_{\varepsilon>0}[C_{\varepsilon
}]_{\infty}}+\sigma_{\cap_{\delta>0}\operatorname*{cl}\left(  \cup
_{0<\varepsilon<\delta}[B_{\varepsilon,\tilde{\rho}_{\varepsilon}}]_{\infty
}\right)  }\nonumber\\
&  =\sigma_{\lbrack\cap_{\varepsilon>0}A_{\varepsilon}]_{\infty}}%
+\sigma_{\lbrack\cap_{\varepsilon>0}C_{\varepsilon}]_{\infty}}+\sigma
_{\cap_{\delta>0}\operatorname*{cl}\left(  \cup_{0<\varepsilon<\delta
}[B_{\varepsilon,\tilde{\rho}_{\varepsilon}}]_{\infty}\right)  }, \label{inn}%
\end{align}
where in the last equality we use\ (\ref{ricone}). Arguing as above, we note
that the three support functions appearing on the right-hand side of
(\ref{inn}) are convex, proper and nonnegative, and are dominated\ by
$\sigma_{\mathrm{N}_{\operatorname*{dom}f}(x)}.$ Thus, they\ are also
continuous on the (nonempty) set $\operatorname*{int}(\operatorname*{dom}%
f)-x$. Applying again Moreau-Rockafellar's Theorem in\ (\ref{inn}) entails
\begin{align*}
\mathrm{N}_{\operatorname*{dom}f}(x)  &  =\partial\sigma_{\mathrm{N}%
_{\operatorname*{dom}f}(x)}(\theta)\\
&  =\partial\sigma_{\left[  \cap_{\varepsilon>0}A_{\varepsilon}\right]
_{\infty}}(\theta)+\partial\sigma_{\left[  \cap_{\varepsilon>0}C_{\varepsilon
}\right]  _{\infty}}(\theta)+\partial\sigma_{\cap_{\delta>0}\operatorname*{cl}%
\left(  \cup_{0<\varepsilon<\delta}[B_{\varepsilon,\tilde{\rho}_{\varepsilon}%
}]_{\infty}\right)  }(\theta)\\
&  =\left[  \cap_{\varepsilon>0}A_{\varepsilon}\right]  _{\infty}+\left[
\cap_{\varepsilon>0}C_{\varepsilon}\right]  _{\infty}+\cap_{\delta
>0}\operatorname*{cl}\left(  \cup_{0<\varepsilon<\delta}[B_{\varepsilon
,\tilde{\rho}_{\varepsilon}}]_{\infty}\right)  .
\end{align*}
Due to the scaling property $\lambda C_{\varepsilon}=C_{\lambda\varepsilon}$
for all $\lambda>0$, it follows that\ $\left[  \cap_{\varepsilon
>0}C_{\varepsilon}\right]  _{\infty}=\{\theta\}\cup(\cap_{\varepsilon
>0}C_{\varepsilon}).$ Consequently, (\ref{f2}) follows from the above equality
together with the fact that\ $[B_{\varepsilon,\tilde{\rho}_{\varepsilon}%
}]_{\infty}=[B_{\varepsilon,\alpha}]_{\infty}$ (see \ref{lc}).
\end{dem}

\medskip Continuity assumptions on the supremum function allow for simpler
characterizations of its subdifferential. For instance, if\ $f$ is continuous
at the nominal point $x\in\operatorname*{dom}f,$ then (\cite{Val69, Vole1};
see, also, \cite{CHLBook})
\[
\partial f(x)=\bigcap_{\varepsilon>0}\overline{\operatorname*{co}}\left(
{\textstyle\bigcup_{t\in T_{\varepsilon}(x)}}
\partial_{\varepsilon}f_{t}(x)\right)  ,
\]
so in this case the normal cone $\mathrm{N}_{\operatorname*{dom}f}(x)$ reduces
to $\{\theta\},$ and only the $\varepsilon$-active data functions matter. In
what follows, we relax this continuity condition, requiring only that $f$ be
continuous at some point of $\operatorname*{dom}f.$ In this case, the sets
$A_{\varepsilon}(x),$ $B_{\varepsilon,\alpha}(x)$ and $C_{\varepsilon}(x),$ as
defined in (\ref{sets}) and (\ref{setsb}), come into play. The
coefficient\ $\rho_{t,\varepsilon}$ is defined in (\ref{vetb}).

\begin{theo}
\label{ola} Let $f_{t}:X\rightarrow\overline{\mathbb{R}},$ $t\in T,$ be convex
and lsc, and assume that\ $f:=\sup_{t\in T}f_{t}$ is continuous at some point.
Then, for every $x\in\operatorname*{dom}f,$%
\begin{align*}
\partial f(x)=\left(  \bigcap_{\varepsilon>0}\overline{\operatorname*{co}%
}(A_{\varepsilon}(x))\right)  +\left(
{\textstyle\bigcap\limits_{\delta>0}}
\operatorname*{cl}\left(
{\textstyle\bigcup\limits_{0<\varepsilon<\delta}}
[\overline{\operatorname*{co}}(B_{\varepsilon,\alpha}(x))]_{\infty}\right)
\right)  \\
+\left(  \{\theta\}\cup%
{\textstyle\bigcap\limits_{\varepsilon>0}}
\overline{\operatorname*{co}}(C_{\varepsilon}(x))\right)  , \label{tayb1}%
\end{align*}
for any\ $\alpha=\alpha(\varepsilon)\in\mathbb{R}_{+}^{\mathcal{P}\backslash
T_{\varepsilon}(x)}$ such that $\alpha_{t}\geq\rho_{t,\varepsilon}$ and both
$\sup_{t}\alpha_{t}$ and $\inf_{t}\alpha_{t}f_{t}(x)$ are finite.
\end{theo}

\begin{dem}
Fix $x\in\operatorname*{dom}f.$ According to \cite{HL08, HLZ08}, we have
\[
\partial f(x)=\mathrm{N}_{\operatorname*{dom}f}(x)+\cap_{\varepsilon
>0}\overline{\operatorname*{co}}\left(  \cup_{t\in T_{\varepsilon}(x)}%
\partial_{\varepsilon}f_{t}(x)\right)  .
\]
Hence, the conclusion follows by applying Theorem\ \ref{conecont}.
\end{dem}

\medskip The following corollary deals with the compact-continuous framework,
assuming that $T$ is compact and that, for every $z\in\operatorname*{dom}f,$
the mappings $t\mapsto f_{t}(z)$ are upper semicontinuous on $T.$ Recall that
$\rho_{t,\varepsilon}$ is defined in (\ref{vetb}), and that $T(x):=\{t\in
T:f_{t}(x)=f(x)\},$ $x\in f^{-1}(\mathbb{R}).$

\begin{cor}
\label{olab}\ Let\ $f_{t}:X\rightarrow\overline{\mathbb{R}},$ $t\in T,$ be
convex and lsc functions, and denote\ $f:=\sup_{t\in T}f_{t}.$ Under the
compact-continuous framework, for every $x\in f^{-1}(\mathbb{R})$\ we have
\[
\partial f(x)=\bigcap_{\varepsilon>0}\overline{\operatorname*{co}}\left(
\tilde{A}_{\varepsilon}(x)+\varepsilon(\tilde{B}_{\varepsilon,\alpha}(x)\cup
C_{\varepsilon}(x)\cup\{\theta\})\right)  ,
\]
where $C_{\varepsilon}(x)$ is defined in (\ref{setsb}),
\[
\tilde{A}_{\varepsilon}(x):=%
{\textstyle\bigcup_{t\in T(x)}}
\partial_{\varepsilon}f_{t}(x)\text{ and }\tilde{B}_{\varepsilon,\alpha}(x):=%
{\textstyle\bigcup_{t\in\mathcal{P}\backslash T(x)}}
~\partial_{\varepsilon}(\alpha_{t}f_{t})(x),
\]
with\ $\alpha\in\mathbb{R}_{+}^{\mathcal{P}\backslash T(x)}$\ such that
$\alpha_{t}\geq\rho_{t,\varepsilon},$ $\sup_{t}\alpha_{t}<+\infty$ and
$\inf_{t}\alpha_{t}f_{t}(x)>-\infty.$ In addition, if $f$ is continuous at
some point, then
\begin{align*}
\partial f(x)=\left(  \bigcap_{\varepsilon>0}\overline{\operatorname*{co}%
}\left(  \tilde{A}_{\varepsilon}(x)\right)  \right)  +\left(
{\textstyle\bigcap\limits_{\delta>0}}
\operatorname*{cl}\left(
{\textstyle\bigcup\limits_{0<\varepsilon<\delta}}
[\overline{\operatorname*{co}}(\tilde{B}_{\varepsilon,\alpha}(x))]_{\infty
}\right)  \right) \\ +\left(  \{\theta\}\cup%
{\textstyle\bigcap\limits_{\varepsilon>0}}
\overline{\operatorname*{co}}(C_{\varepsilon}(x))\right)  .
\end{align*}

\end{cor}

\begin{dem}
The first conclusion follows similarly to the proof of Theorem \ref{thmsub},
based on the characterization (\cite{HL08, HLZ08})
\[
\partial f(x)=\cap_{\varepsilon>0,L\in\mathcal{F}(x)}\overline
{\operatorname*{co}}\left(  \tilde{A}_{\varepsilon}(x)+\mathrm{N}%
_{L\cap\operatorname*{dom}f}(x)\right)  .
\]
Under the additional continuity assumption, we have
\[
\partial f(x)=\mathrm{N}_{\operatorname*{dom}f}(x)+\cap_{\varepsilon
>0}\overline{\operatorname*{co}}\left(  \tilde{A}_{\varepsilon}(x)\right)  ,
\]
and the second conclusion follows as in the proof of Theorem \ref{ola}.
\end{dem}

We apply\ Theorem\ \ref{ola} in\ the affine case.

\begin{exam}
\label{exammcont}Given a nonempty set $T,$ we consider the function
\[
f:=\sup_{t\in T}\left\{  \left\langle a_{t},\cdot\right\rangle -b_{t}\right\}
,\text{ }a_{t}\in X^{\ast},\text{ }b_{t}\in\mathbb{R}.
\]
We assume that $f$ is continuous at some point. Then, by Theorem
\ref{ola},\ for all $x\in f^{-1}(0)$ we have
\[
\partial f(x)=\left(  \cap_{\varepsilon>0}\overline{\operatorname*{co}%
}(\left\{  a_{t}:\left\langle a_{t},x\right\rangle -b_{t}\geq-\varepsilon
\right\}  )\right)  +\cap_{\delta>0}\operatorname*{cl}\left(  \cup
_{0<\varepsilon<\delta}[\overline{\operatorname*{co}}(B_{\varepsilon
}(x))]_{\infty}\right)  .
\]
where $B_{\varepsilon}(x):=\left\{  \varepsilon(2b_{t}-2\left\langle
a_{t},x\right\rangle -\varepsilon)^{-1}a_{t}:\left\langle a_{t},x\right\rangle
-b_{t}<-\varepsilon\right\}  .$
\end{exam}

\section{Concluding remarks}

This paper establishes a comprehensive and explicit framework for the
subdifferential of pointwise suprema of convex functions, expressed directly
in terms of the underlying data functions. By treating $\varepsilon$-active
and non-$\varepsilon$-active functions symmetrically, the proposed
characterization unifies and extends previous results and avoids cumbersome
geometric constructions. In particular, it yields sharp Karush--Kuhn--Tucker
and Fritz--John optimality conditions for infinite convex optimization,
clearly distinguishing the roles of $\varepsilon$-active and non-$\varepsilon
$-active constraints.


\bibliographystyle{siamplain}

\end{document}